\documentclass[
11pt,%
tightenlines,%
notitlepage,%
twoside,%
onecolumn,%
nobibnotes,%
nofootinbib,%
superscriptaddress,%
noshowpacs,%
noshowydk,%
showdoi,%
centertags,%
maik]%
{revtex4n}

\input jpack_e
\setcaptionmargin{5mm}
\onelinecaptionstrue

\columntovsizetrue
\allowdisplaybreaks

\begin{document}

\newtheorem{teo}{Theorem}
\newtheorem*{teon}{Theorem}
\newtheorem{lem}{Lemma}
\newtheorem*{lemn}{Lemma}
\newtheorem{prp}{Proposition}
\newtheorem*{prpn}{Proposition}
\newtheorem{ass}{Assertion}
\newtheorem*{assn}{Assertion}
\newtheorem{assum}{Assumption}
\newtheorem*{assumn}{Assumption}
\newtheorem{stat}{Statement}
\newtheorem*{statn}{Statement}
\newtheorem{cor}{Corollary}
\newtheorem*{corn}{Corollary}
\newtheorem{hyp}{Hypothesis}
\newtheorem*{hypn}{Hypothesis}
\newtheorem{con}{Conjecture}
\newtheorem*{conn}{Conjecture}
\newtheorem{dfn}{Definition}
\newtheorem*{dfnn}{Definition}
\newtheorem{problem}{Problem}
\newtheorem*{problemn}{Problem}
\newtheorem{notat}{Notation}
\newtheorem*{notatn}{Notation}
\newtheorem{quest}{Question}
\newtheorem*{questn}{Question}

\theorembodyfont{\rm}
\newtheorem{rem}{Remark}
\newtheorem*{remn}{Remark}
\newtheorem{exa}{Example}
\newtheorem*{exan}{Example}
\newtheorem{cas}{Case}
\newtheorem*{casn}{Case}
\newtheorem{claim}{Claim}
\newtheorem*{claimn}{Claim}
\newtheorem{com}{Comment}
\newtheorem*{comn}{Comment}

\theoremheaderfont{\it}
\theorembodyfont{\rm}

\newtheorem{proof}{Proof}
\newtheorem*{proofn}{Proof}

\selectlanguage{english}%%
\Rubrika{\relax}
\CRubrika{\relax}
\SubRubrika{\relax}
\CSubRubrika{\relax}
%

% Current Number and Volume of Journal
\def\JournalNumber{}
\def\JournalVolume{}
%
%
% RegDyn_01_07.tex 20.12.2006
%
\nameVolumeRus{}%    => \ContentsHeadLineRusA
\CnameVolumeRus{}%   => \ContentsHeadLineRusA
\nameIssueRus{\No}%  => \ContentsHeadLineRusA
\CnameIssueRus{}%    => \ContentsHeadLineRusA
\namePartRus{}%      => Не используется
\namePagesRus{}%     => \ArticleDataHeadRusA
\nameYearShortRus{}% => Не используется
\JournalNameRus{}%   => \ArticleDataHeadRusA
\TranslitJournalNameRus{}% => \ArticleDataHeadB; \@ArticleDataFootRus
\JournalName{}% \ContentsHeadLineC; \ArticleDataHeadA; \@oddfoot; \@evenfoot
\JournalISSNCode{0000-0000}% Указывать обязательно, используется при расчете DOI => \ContentsHeadLineC
\IssuePrice{}% => Не используется
\TransYearOfIssue{2018}% => \ContentsHeadLineA; \ArticleDataHeadA; \@oddfoot; \@evenfoot
%\TransCopyrightYear{}%
\OrigYearOfIssue{}% => \ContentsHeadLineRusA; \ArticleDataHeadRusA; \ArticleDataHeadC; \@ArticleDataFootRus
\OrigCopyrightYear{}%
\OrigIssueNo{\JournalNumber}% \ArticleDataHeadRusA;
\OrigVolumeNo{\JournalVolume}% => \ContentsHeadLineRusA; \ArticleDataHeadRusA; \ArticleDataHeadA
\TransVolumeNo{\JournalVolume}% => \ContentsHeadLineA; \ArticleDataHeadA; \@oddfoot; \@evenfoot
\TransIssueNo{\JournalNumber}% => \ContentsHeadLineA; \ArticleDataHeadA; \ArticleDataHeadRusA; \@oddfoot; \@evenfoot
\TransPartNo{}% => Не используется
\SHORTjournalPREFIX{} % Короткое имя журнала
\LONGjournalPREFIX{} % Длинное имя журнала
\BatFileName{call make_ps.bat} % имя командного файла для производства ps-файлов
\BatSwitch{3} % ps-файлы по номерам (не менять)
%%%%%%%%%%%%%%%%%%%%%%%%%%%%%%%%%%%%%
%\IssueName{}% Наименование выпуска
\SupplementNumber{}% Номер yes/no --- SGML
\PublicationSerialNumberInYear{0}% Порядковый номер издания в году --- используется при расчете DOI
\PublicationSerialNumberInVolume{0}% Порядковый номер издания в томе
%\ConditionalIssueDate{"year","month","day","name","type"}% Условная дата выпуска  --- SGML
%\PagePrefix{}% --- SGML
%\JournalISSNonlineCode{}%  --- SGML
%\JournalISSNCodeRus{}% --- SGML
%\JournalISSNonlineCodeRus{}% --- SGML
%%%%%%%%%%%%%%%%%%%%%%%%%%%%%%%%%%%%%
%\VolumeName{}% --- SGML
%\IssnoName{none}% --- SGML
%\PartnoName{}% --- SGML
%\FpageNamepp{}% --- SGML
%\FpageNnamep{}% --- SGML
%\FpagePrefix{}% --- SGML
%\LpageNnamepp{}% --- SGML
%\LpageNamep{}% --- SGML
%\LpagePrefix{}% --- SGML
%\VolumePageNumbering{}% yes/no --- SGML
%\JournalPubID{}% --- SGML
%\FirstJournalPageNumber{}% --- SGML
%\LastJournalPageNumber{}% --- SGML
%%%%%%%%%%%%%%%%%%%%%%%%%%%%%%%%%%%%%
\makeatletter
\def\MAIKlogo{}%  => \ContentsHeadLineC Pleiades~Publishing,~Ltd.
\def\maikpraefix{}% MAIK; Pleiades~Publishing,~Inc.
%\def\maikpraefix{10.0000/S}% MAIK; Pleiades~Publishing,~Inc.
% Вторая строка оглавления
%\edef\@ContentsHeadLineB{Simultaneous English language translation of the journal is available from \noexpand\MAIKlogo}
% Начало третьей строки оглавления
%\def\Distributed{Distributed worldwide by Springer. }% => \ContentsHeadLineC
%
%\def\ArticlePages#1{\relax}
%\@ifxundefined\CONT@sw{\@booleantrue\CONT@sw}{}%
%\@booleantrue\showPACS@sw%
%\@booleanfalse\showPACS@sw
\@booleantrue\showKEYS@sw %
%\@booleantrue\noOrigJournalVersion@sw% Не русской версии журнала --- отсутстует
%\@booleantrue\noOrigVolumeNo@sw% Не выводится том русской версии в колонтитулах
%\@booleanfalse\noTransVolumeNo@sw% Выводится том английской версии в колонтитулах и оглавлении
%\makeatother
%
 
 %
\input maikdoi %

\beginpaper

%*************************************************************

\input engnames%% For article in English
\titlerunning{Quasiperiodic orbits in Siegel disks/balls and the Babylonian problem}%Short title in headings
\authorrunning{Y. Saiki and J.~A. Yorke}%Author in headings -
\toctitle{Quasiperiodic orbits in Siegel disks/balls and the Babylonian problem}%Title in contents
\tocauthor{Y. Saiki and J.~A. Yorke}%Author in contents
\title{Quasiperiodic orbits in Siegel disks/balls and the Babylonian problem}%Full Title article
\firstaffiliation{%% Affiliation of organization
}%
\articleinenglish %%
\PublishedInRussianNo
\author{\firstname{Yoshitaka}~\surname{Saiki}}%
\email[E-mail: ]{yoshi.saiki@r.hit-u.ac.jp}
\affiliation{%% Affiliation1
Graduate School of Business Administration, Hitotsubashi University\\
2-1 Naka, Kunitachi, Tokyo 186-8601, Japan}%
\affiliation{%% Affiliation1
JST PRESTO\\
4-1-8 Honcho, Kawaguchi-shi, Saitama 332-0012, Japan}
\affiliation{%% Affiliation2
University of Maryland\\
College Park, MD 20742, USA}%
\author{\firstname{James A.}~\surname{Yorke}}%
\email[E-mail: ]{yorke@umd.edu}
\affiliation{%% Affiliation2
University of Maryland\\
College Park, MD 20742, USA}%
%\noaffiliation
\begin{abstract}%%
We investigate numerically complex dynamical systems where a fixed point is surrounded by a disk or ball of quasiperiodic orbits,
where there is a change of variables (or conjugacy) that converts the system into a linear map. We compute this  ``linearization'' (or conjugacy) from knowledge of a single quasiperiodic trajectory. 
In our computations of rotation rates of the almost periodic orbits and Fourier coefficients of the conjugacy, we only use knowledge of a trajectory, and  
we do not assume knowledge of the explicit form of a dynamical system.
This problem is called the Babylonian Problem: 
determining the characteristics of a quasiperiodic set from a trajectory. Our computation of rotation rates and Fourier coefficients depends on the 
very high speed of our computational method  ``the weighted Birkhoff average''.
%In the abstract of the paper you should describe the system under
%consideration, the methods and approaches used and the results obtained. As a
%rule, the volume of the abstract should be at least 100 words.
\end{abstract}
%MSC2010 -
%\received{September 24, 2018}%please insert here date of submitting
%\accepted{October 30, 2018}%
\maketitle

%\textmakefnmark{0}{)}%

%Here You may place your own definition
%\usepackage{comment}
\newcommand{\blue}{\textcolor{blue}}
\newcommand{\blu}{\textcolor{blue}}
\newcommand{\magenta}{\textcolor{magenta}}
\newcommand{\red}{\textcolor{red}}
\newcommand{\black}{\textcolor{black}}
\newcommand{\green}{\textcolor{green}}
\newcommand{\cya}{\textcolor{cyan}}
\newcommand{\cyan}{\textcolor{cyan}}
\newcommand{\Q}{\mbox{WB}}
\newcommand{\B}{\mbox{B}}
\newcommand{\WB}{\mbox{WB}}
\newcommand{\Trho}{\rho}
\newcommand{\torus}{{\mathbb{T}^d}}
\newcommand{\TTT}{{\mathbb{T}}}
\newcommand{\Torus}{{\mathbb{T}^2}}
\newcommand{\QC}{\mbox{zeros}}
\newcommand{\fd}{\mbox{fd}}%for \int\fd\mu
\newcommand{\fk}{e^{\Iota 2\pi k\cdot\theta}}%{f_k(\theta)}
\newcommand{\R}{{\mathbb{R}}}
\newcommand{\mm}{\mathbb}
\newcommand{\photo}{png}
\newcommand{\gd}{\bar\Delta^*}
\newcommand{\mR}{\mathbb{R}}
\newcommand{\allblue}{\color{black}{}}
\newcommand{\allred}{\color{red}{}}
\newcommand{\allblack}{\color{black}{}}
\newcommand{\allmagenta}{\color{magenta}{}}
\newcommand{\Path}{./figures/}
\newcommand{\Xx}{\ref{eq:plusrho}}
\newcommand{\D}{\mathcal{D}}
\newcommand{\BF}{\mathbf}
\renewcommand{\thefootnote}{\fnsymbol{footnote}}
\newcommand\keywordsname{Key words}
\newcommand\AMSname{AMS subject classifications}
%.................................................................................%

%\Keywords{quasi-periodic orbit; rotation rate; Siegel disk; Siegel ball}

\renewcommand*\contentsname{Contents}
%\tableofcontents
%\allred
%rotation rate? or rotation rate?
%\allblack

\section{Introduction}\label{s:introduction}
Let $F:\mathbb{C}^d\to\mathbb{C}^d$ be a complex analytic map in $d$ complex dimensions.
There are situations where there is an open set of points $z_0$ for which the trajectory $z_n = F^n(z_0)$ for $n\in \mm N$ is quasiperiodic. 
%\thefootnote{test}

{\bf The Babylonian Problem.}
Our goal is to obtain information about the dynamics only from knowledge of a single quasiperiodic  trajectory $(z_n)$, using no additional information about $F$. Our emphasis is on numerical examples for $d=1$ and $2$.  We named the task of determining information about the dynamics -- such as determining the rotation rates 
-- from a trajectory the {\bf Babylonian Problem}, named after the Babylonians who about 2500 years ago determined the three rotation rates of the orbit of the moon; see \cite{Goldstein}. The moon's orbit can be viewed as approximately as quasiperiodic on a three-dimensional torus. (That ignores other frequencies as caused by the eccentricity of the Earth's orbit around the sun and the influence of other planets, especially Jupiter.) The data they used was the trajectory of the moon viewed against the fixed stars.
%REFERENCE OUR BABYLONIAN PAPER HERE.  
See \cite{DCDS} for the details of our method on the Babylonian problem.
%Section\ref{rotation-vector} for a brief summary and 
%\cite{DCDS} for the details.
\allblack

We study quasiperiodic orbits in a one-dimensional complex dynamical system in Section~\ref{s:siegel-disk}, 
and those in two-dimensional complex H\'enon map in Section~\ref{s:siegel-ball}, but we introduce them here.

\allblue
{\bf Quasiperiodicity defined.}\label{ss:quasiperiodicity-def} 
Let $\torus$ be a $d$-dimensional torus. We will frequently use the coordinate representation 
$\theta \in\torus
:= [0,1]^d \bmod1$ where $\bmod~1$ is applied to each coordinate.
For $\rho\in\torus$, we will say $\rho=(\rho_1,\cdots,\rho_d)$ is {\bf irrational} when all coordinates $\rho_k$ are irrational and {\bf rationally independent}, 
i.e., if whenever for $k=1,\cdots,d,$ $a_k$ are integers  for which 
$a_1\rho_1+\cdots+a_d\rho_d = 0$, then $a_k=0$ for all $k=1,\cdots,d$. 

The simplest example of a quasiperiodic map is
\begin{equation}\label{eq:plusrho}
\theta_{n+1} := F(\theta_n) = \theta_n + \rho\bmod1
\end{equation}
where $\rho$ is irrational. We will refer to $\rho$ as a 
{\bf rotation vector}. 

A map is {\bf $d$-dimensionally quasiperiodic} when there  exists a smooth choice of coordinates for
the torus such that the dynamics on the orbit is 
given by the map Eq.~\ref{eq:plusrho}, where $\rho$ is irrational.

\allblack

%
%%\ifincludeXX
\begin{figure}
  \centering{
    \includegraphics[width = 0.45\textwidth,height=.35\textwidth]{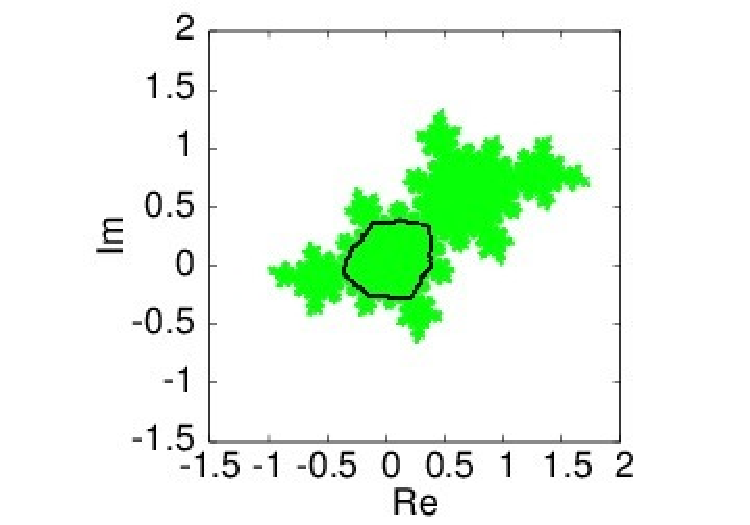}
    \includegraphics[width = 0.45\textwidth,height=.35\textwidth]{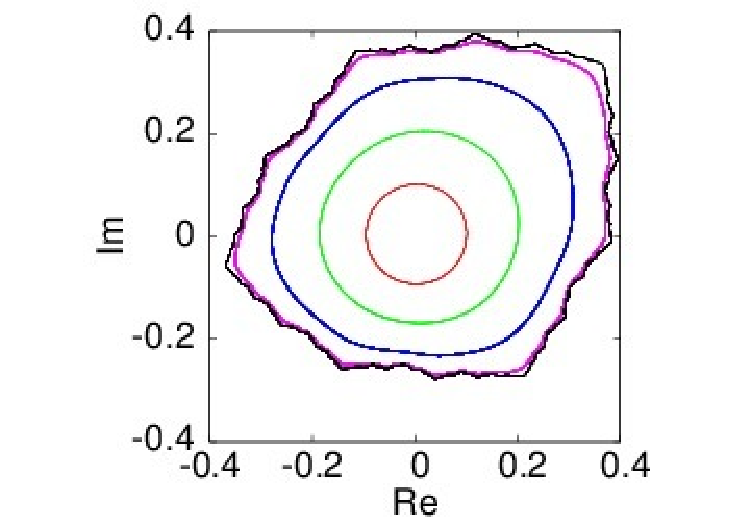}
    \caption{\textbf{Quasiperiodic orbits in the Siegel disk of Example 1 in Eq.~\ref{eq:1d}.} Left: A bounded curve (black) is plotted on top of a set of points (green) whose orbits stay bounded. The fractal boundary of the set is a Julia set. 
Right:    Five bounded curves are shown, and four curves except for the outermost curve are focused on in Figs.~\ref{fig:fourier-siegel-converted}, ~\ref{fig:fourier-siegel} and ~\ref{fig:convergence-siegel}. 
The curves are generated from initial points $z_0$=$0.1$ (red), $0.2$ (green), $0.3$ (blue), $0.37$ (magenta) in $\mathbb{C}$, 
and the outermost curve (black) is 
derived from the Fourier coefficients
Eq.~\ref{eq:rfourier} 
%Eq.~\ref{eq:h-power-series} 
of the curve from the initial point $z_0=0.37$, the image under $h$ of the circle at $r=0.999$ in the unit disk $\D$. The outermost curve in the Right panel is the same as the one plotted in the Left panel.
    }\label{fig:quasi-siegel}
}
\end{figure}
%\fi

%\ifincludeXX
\begin{figure}
  \centering
  \includegraphics[width = .45\textwidth,height=.35\textwidth]{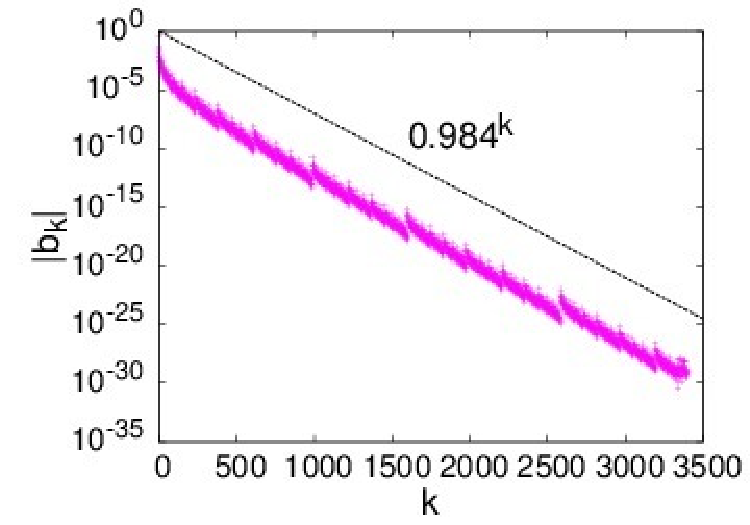}
  \includegraphics[width = .45\textwidth,height=.35\textwidth]{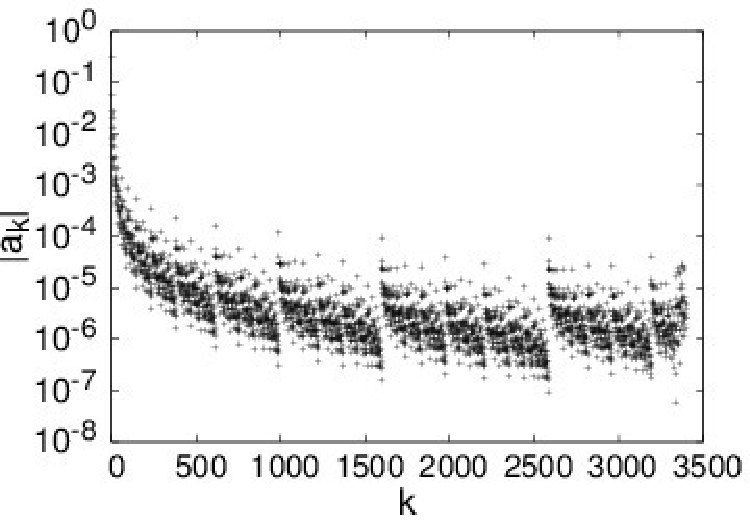}
  \caption{
    \textbf{Fourier and Power series coefficients.}
      Left: Fourier coefficients $b_k$ for a curve starting from 
      $z_0=0.37$ shown in Fig.~\ref{fig:quasi-siegel}~(Right).
    Its slope is $R_0\sim 0.984$.
    Right: Power series coefficients $a_k:=R_0^{-k}b_k$ of the conjugacy $h$ in Eq.~\ref{eq:conjugacy}. See also Eqs.~\ref{eq:R0} and \ref{eq:h-power-series}.
    We used $N=4\cdot10^6$ iterates $(z_n)$ for the computation of Fourier coefficients.
}
\label{fig:fourier-siegel-converted}
\end{figure}
%\fi
%

\subsection {Our primary example in $\mathbb{C}$} Define
\begin{equation}
  F(z)
  =F_1(z)
  :=z^2+e^{2\pi i \rho}z \mbox{ where }\rho=(\sqrt{5}-1)/2,  
\label{eq:1d}
\end{equation}
so ${dF(0)}/{dz}=e^{2\pi i \rho}$. 
Fig.~\ref{fig:quasi-siegel} shows several trajectories, each of which is dense on a closed curve. 
\allblue
From the above definition, for $d=1$, saying that 
a trajectory of Eq.~\ref{eq:1d} is quasiperiodic is equivalent to saying that there is a conjugacy
$H:S^1\to\mathbb{C}$ that maps the trajectory $(e^{i2\pi n \rho})_{n\in\mm N}$ smoothly onto $(z_n)$ for irrational $\rho$.
More specifically, again interpreting $S^1$ as the unit circle in $\mathbb{C}$,
\allblack
\begin{equation}
F(H(e^{i2\pi\theta}))= H(e^{i2\pi(\theta+\rho)}) \mbox{ for each }e^{i2\pi\theta}\in S^1.
\end{equation}
%\allred
See De la Llave 2018~\cite{delaLlave2018}, Siegel and Moser 1995 (\cite{Siegel-Moser} P.185) for discussions of the conjugacy.
%Perhaps this is redundant, but the original text also suggested Siegel and Carleson.
%\allblack

For the Babylonian problem we first determine the rotation rate $\rho$ from $(z_n)$. 
Our paper \cite{DCDS} solved the problem of determining $\rho$ from a quasiperiodic trajectory; methods for some cases had previously been established, but our method works in all cases with smooth dynamics including higher dimensional cases. See
\cite{DCDS, Luque-Villanueva2009} 
and references therein.

Secondly we can determine $H$.
Using our group's ``weighted Birkhoff'' averaging methods described in Sec.~\ref{ss:WB}, we compute the Fourier series coefficients 
 from $(z_n)$
 for the conjugacy $H:S^1\to \mathbb{C}$, yielding 
\begin{equation}
H(e^{i2\pi\theta}) =\sum_{k=0}^{+\infty} b_k 
\sigma_k(\theta) \mbox{ where }\sigma_k(\theta):=
e^{i 2\pi k \theta}.\label{eq:H}
\end{equation}
%In Sect.~\ref{ss:WB} we outline our group's ``weighted Birkhoff'' method for computing the rotation rates $\rho$ and the coefficients $b_k$ in Eq.~\ref{eq:H}. 
Our method is useful because of its very fast rate of convergence, especially useful when computing thousands of $b_k$ values.
We can compute (or rather estimate) $\rho$ and each $b_k$ from $N$ trajectory points of $(z_n)_{n=1}^N$.
Write $\rho^{(N)}$ and $b_k^{(N)}$ for our approximations obtained using the $N$ terms  $(z_n)_{n=0}^N$.
The computational errors are $|\rho^{(N)}-\rho| \mbox{ and } |b_k^{(N)}-b_k|$.
We have proved that when $\rho$ satisfies a Diophantine condition, Ineq.~\ref{ineq:diophantine},
we get ``super-fast'' (or ``faster than polynomial'') convergence, as follows.
\begin{equation}
\mbox{ For each }m\in\mm Z,~
|\rho^{(N)}-\rho|\cdot N^m\to 0 \mbox{ and } |b_k^{(N)}-b_k|\cdot N^m\to 0 \mbox{ as }N\to\infty   \label{super}
\end{equation}
for all $k\in\mm Z$. See Theorem~\ref{thm:A1} for a more general statement.

\subsection{Siegel Disk and the conjugacy.} 
Let  $U$ be an open subset of $\mathbb{C}$ containing $0$, and $F: U \to \mathbb{C}$ be holomorphic with $F(0)=0$ and with 
${dF}/{dz}(0) = e^{i2\pi \rho}$ for some real irrational $\rho$. 
%\allred
%satisfying a certain condition discussed later.
%\allblack 
Siegel showed that for almost all $\rho$, $F$ is {\bf linearizable} in some neighborhood $\mathcal{U}$ of $0$ in the sense that
$F$ is conjugate to the linear map 
$\phi(z) = e^{i2\pi \rho}z$. 
That in turn means there is a holomorphic map $h$ from the {\bf open unit disk} $\D$ onto $\mathcal{U}$ for which 
\begin{equation}
F\circ h = h\circ \phi.\label{eq:conjugacy}
\end{equation} 
Note that $F(\D)=U$. 
The maximal such neighborhood $\mathcal{U}$ is called the {\bf Siegel disk} of $0$.
See Siegel \cite{Siegel42}. Also see 
\cite{Luque-Villanueva2009,Llave-Petrov2008, Milnor}
%\cite{Llave-Petrov2008, Luque-Villanueva2009,Ushiki2016}
and references therein for results on the Siegel disk including numerical computations.

{\bf Power series for the conjugacy 
%$h$
$\mathbf{h}$
 from a Fourier series.}
For almost any $\rho$, from one trajectory  $z_{n+1}=F(z_n)$ for $n\ge0$ on $\mathcal{U}$, we can numerically determine the 
power series of the conjugacy map $h$ with high numerical precision,  (without any other knowledge of 
$F$). 
Let $R_0$ be the infimum of $R>0$ for which the sequence
$(b_k/R^k)_{k\ge 0}$  is bounded. 
Of course this depends on the curve that is chosen and we have chosen one with 3400 coefficients that have absolute value above $10^{-30}.$ It appears that $|b_k|<10^{-28}$ for $k>3400$.
We identify its $R_0$ value as $\sim 0.984$. 
Since we only have the finite number of $b_k$, 
we cannot determine $R_0$ value with more than about 3 digits of precision.
Define 
\begin{equation}
a_k := b_k/R_0^k.
\label{eq:R0}
\end{equation}
Then $h$ is an analytic function when defined as
\begin{equation}
h(z):= \sum_{k=0}^{+\infty} a_k z^k \mbox{ for }|z|<1.\label{eq:h-power-series}
\end{equation}
%
%
%\ifincludeXX
\begin{figure}
  \centering
  \includegraphics[width = .32\textwidth,height=.35\textwidth]{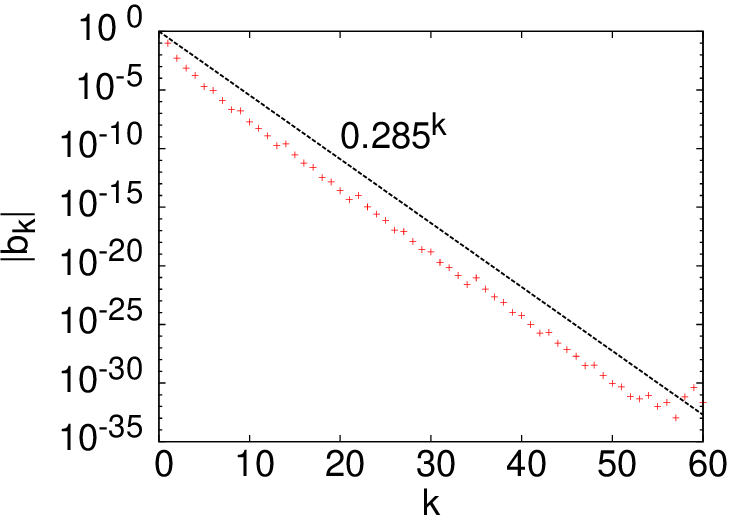}
  \includegraphics[width = .32\textwidth,height=.35\textwidth]{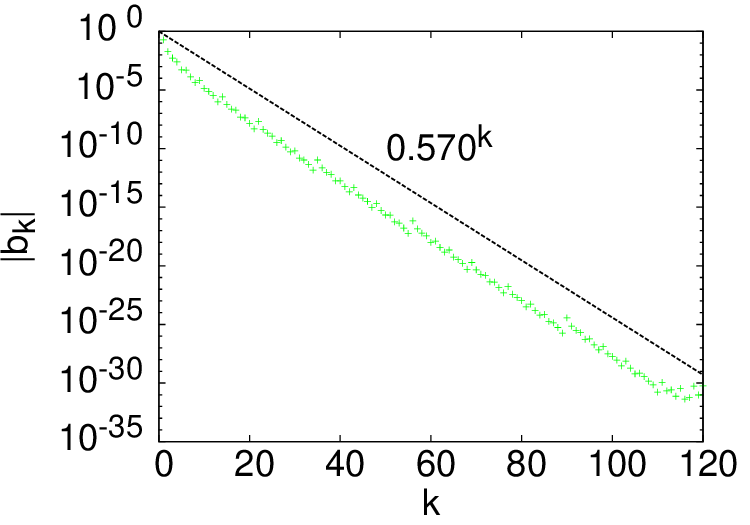}
    \includegraphics[width = .32\textwidth,height=.35\textwidth]{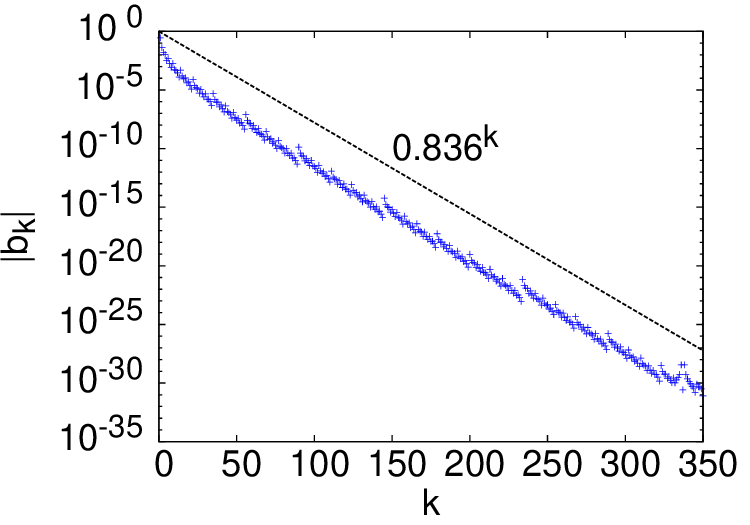}
    \caption{
      \textbf{Magnitude of the Fourier coefficients for three quasiperiodic orbits (Fig.~\ref{fig:quasi-siegel}~(Right)) for Eq.~\ref{eq:1d}.} The initial points $z_0 \in \mathbb{C}$ are $z_0=0.1$~(Left), $0.2$~(Middle), $0.3$~(Right).
     Results for $z_0=0.37$ are shown in 
    Fig.~\ref{fig:fourier-siegel-converted}.
      The slope of each curve corresponds to the radius $r$ of the curve's pre-image under $h$ within $\D$.  
      Note that all three curves show evidence of irregular behavior in the region where $|b_k|\sim 10^{-30}$ reflecting the limits of quadruple precision.
For the computations of the Fourier coefficients  the $N=10^5$ iterates for $0.1$ and $0.2$ curves and the $10^6$ iterates for $0.3$ curve are used.
      }
  \label{fig:fourier-siegel}
\end{figure}
%\fi
%
%
%\ifincludeXX
\begin{figure}
    \centering{       
            \includegraphics[width = 0.5\textwidth,height=0.4\textwidth]{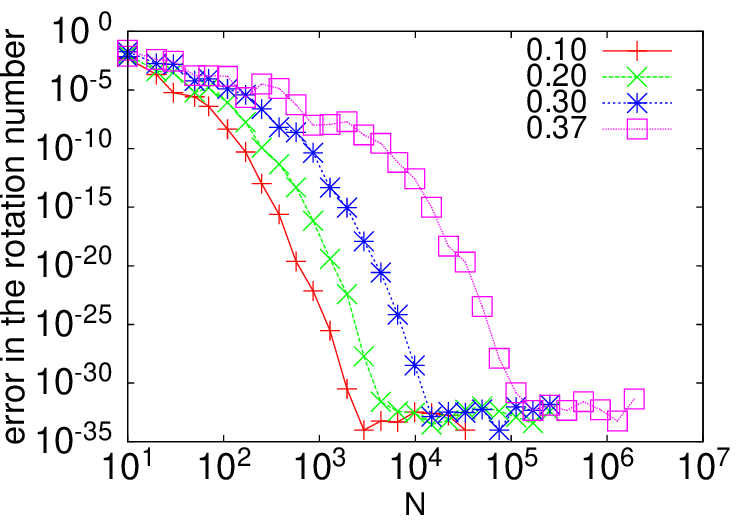}	
      \caption{
        {\bf Convergence of the rotation rate of each of the inner four trajectories in Fig.~\ref{fig:quasi-siegel}~(Right) to the rigorous value ($\rho=(\sqrt{5}-1)/2$)}. The error in the rotation rate $err(N)=|\rho_N -\rho|$ is plotted with respect to the number $N$ of iterates of the trajectory used by the weighted Birkhoff average 
%        $\Q_N^{[1]}$ (left) and 
($\Q_N^{[2]}$ in Sec.~\ref{ss:WB}).
% (right), respectively.}
      }
      \label{fig:convergence-siegel}
}
\end{figure}
%\fi

%
\noindent
%With only 3400 coefficients $b_k$ (see Fig.~\ref{fig:fourier-siegel-converted}~(left)) we can only estimate $R_0$ with an expected three-digit precision.
For $|z|<1$ the above terms in the sum converges exponentially fast to $0$ and 
for $|z|>1$ the terms diverge.
Let $C_r$ denote the circle of radius $r$ centered at $0$.
For $z\in C_{r}$ for $0<r<1$ we can write 
$z=re^{i2\pi\theta}$ in polar coordinates and then 
\begin{equation}
h(re^{i2\pi\theta}) =\sum_{k=0}^{+\infty} a_k r^k 
\sigma_k(\theta) = \sum_{k=0}^{+\infty} b_k (r/R_0)^k
\sigma_k(\theta) .\label{eq:rfourier}
\end{equation}
Hence $h:\D\to\mathcal{U}$ is the desired linearization, a conjugacy of the unit disk to the Siegel Disk. 
We can restrict $h$ to circles $C_r$ of different radii $r<1$ and get good approximations of the corresponding curves in $\mathcal{U}$.
The outermost curve for $r=0.999$ in Fig.~\ref{fig:quasi-siegel} (Right) is computed from Eq.~\ref{eq:rfourier} using $(b_k)_{k=0,\cdots, 3400}$ for a quasiperiodic orbit from $z_0=0.37$, and the estimated value of $R_0\sim 0.984$.
%{\bf Example in $\mathbb{C}$: Example 1.} Our primary one-dimensional example is Eq.~\ref{eq:1d}.

{\bf Diophantine Condition.}
In order to discuss more general values of $\rho$ and to explain which $\rho$ values the map in Eq.~\ref{eq:1d} has nice behavior, Siegel assumed $\rho$ satisfied a Diophantine condition.
For a fixed real number $\kappa \ge 2$ we say that an irrational $\rho\in\R$ satisfies a 
{\bf Diophantine condition of order} $\kappa$ if there exists some $\varepsilon>0$ so that 
\begin{equation}
\left|\rho-\frac{p}{q}\right|> \frac{\varepsilon}{q^\kappa}
\label{ineq:diophantine}
\end{equation}
for every rational rate $p/q$.
%where $p, q \in \mathbb{Z}$ and $q\ne0$. 
Let $\lambda=e^{2\pi i \rho}$. If $\rho$ satisfies Ineq.~\ref{ineq:diophantine}, it follows that
\begin{equation}
|\lambda^q-1|>\varepsilon^{\prime}/q^{\kappa -1}
\end{equation}
for some $\varepsilon^{\prime}>0$ for all positive integers $q$.
Let $D_\kappa$ be the set of all (irrational) numbers $\rho$ which satisfy this condition.
Note that $D_\kappa \subset D_\eta$ whenever $\kappa < \eta$.
\allblue
The set of Diophantine numbers is the union of the $D_\kappa$.
\allblack
 In 1942 Siegel~\cite{Siegel42} proved that if the Diophantine condition holds, the Siegel disk exists for a class of dynamical systems including Eq.~\ref{eq:1d}. 
See Brjuno \cite{Brjuno, Brjuno2} for a more general condition, which is proved to be sharp by Yoccoz~\cite{Yoccoz}.
Note that in our examples we choose parameters which satisfy the Diophantine condition (and of course, Brjuno condition as well).

\subsection{Example in  $\mathbf{\mathbb{C}^2}$: The complex H\'enon map} 
Our main example for $d=2$ follows Ushiki \cite{Ushiki0}. 
%Following Ushiki \cite{Ushiki2016,Ushiki2}, 
We study the Siegel ball of the complex H\'enon map 
$
  F_2: \mathbb{C}^2\to \mathbb{C}^2
$
%\end{equation}
depending on fixed parameters $\alpha, \beta\in\mathbb{C}$, where $F_2$ is given by 
\begin{equation}\label{henon}
F=F_2
    \begin{pmatrix}
      x\\
      y
    \end{pmatrix}
    :=
    \begin{pmatrix}
      y\\
      \beta(y^2+\alpha)-\beta^2 x
    \end{pmatrix}.
\end{equation}
We will specify $\alpha$ and $\beta$ below so that the Jacobian of one of the fixed points has eigenvalues of the form 
$e^{i(\theta + \phi)}$ and $e^{i(\theta - \phi)}$, 
both on the unit circle in $\mathbb{C}$. When  $\theta/ 2\pi$ and $\phi/ 2\pi$ and $\theta/\phi$ are irrational, trajectories of the linearized system about that fixed point are quasiperiodic. 
When $\theta/ 2\pi$ and $\phi/ 2\pi$ and $\theta/\phi$ are Diophantine, 
there is a conjugacy between the linearized system and the dynamics in a neighborhood (called the Siegel ball) of that fixed point.
We could in principle use our weighted Birkhoff method (See Sec.~\ref{ss:WB}) to compute Fourier coefficients of the conjugacy, though this two-dimensional case requires many more coefficients than the one-dimensional case above to obtain the same numerical accuracy.
Instead we examine the curious aspect of determining the rotation vectors.
In $\mathbb{C}^2$ a 2D-quasiperiodic orbit has a conjugacy
$H:S^1\times\ S^1\to \mathbb{C}^2$, which can be written as
$$
H(e^{i\theta_1 2\pi},e^{i\theta_2 2\pi})\mbox{ for } (e^{i\theta_1 2\pi},e^{i\theta_2 2\pi})\in S^1\times S^1.
$$
If $z_0:=H(1,1)$, then there exist $ \sigma_1, \sigma_2$ such that the trajectory is
$$
F^n(z_0) = H(e^{n~i \sigma_1 2\pi},e^{n~i \sigma_2 2\pi}).
$$
The curious aspect of this is that there is a dense set of 
$$
(\sigma_1,\sigma_2)\in [0,1)\times [0,1) \bmod1 
$$
(where $\bmod1 $ applies to each coordinate) 
for which there is such a conjugacy map $H$ that yields the same trajectory $(F^n(z_0))$. 
That is because there are countably many linear choices of coordinates on $S^1\times S^1$.
%\allblue
We resolve this ambiguity by procedures analogous to the Babylonians'. They observed a projection of the moons orbit onto the two-dimensional celestial sphere. We also require a projection.
%\allblack
See Sec.~\ref{ss:phi} below.
% we summarize our results in Section~\ref{s:summary}. 
Furthermore we examine the two-dimensional case with domain in $\mathbb{C}^2$ with analogous Siegel results. 
%\allblue

{\bf Rotation rates for quasiperiodic orbits on }${\BF \torus}$
{\bf are not well defined.} Every quasiperiodic orbit is conjugate
to the map in Eq.~\ref{eq:plusrho}; i.e.,
$$
F(\theta)=\theta+\rho \bmod1
$$
so we might be encouraged to call the vector $\rho\in\torus$ the rotation vector of the map; however, it is also conjugate to other maps.  The matrix $A$ is a {\bf unimodular} transformation if it is an invertible $d\times d$ matrix with integer coefficients. Such a matrix can be viewed as a conjugacy on $\torus$ by writing  
$\bar\theta := A\theta$.
In these coordinates, Eq.~\ref{eq:plusrho} becomes
\begin{equation}\label{eq:rhobar}
\bar\theta_{n+1} = \bar\theta_n + A\rho\bmod1.
\end{equation}
Note that $\rho$ is irrational 
if and only if $A\rho$ is. Hence the irrationality of $\rho$ 
is well defined. 

\allblue
If $\rho$ is a rotation vector, then so is $A\rho$ for the same process -- for every unimodular $A$.
However, as we discuss in \cite{DCDS}, for a given irrational $\rho$
the set  $\{A\rho: A\mbox{ is unimodular}\}$ is dense in $\torus$.
We might want to know the vector $\rho$ with 30-digit precision. But every vector in $\torus$ is a 30-digit approximation of $A\rho$
for an appropriate choice of coordinate matrix $A$.
And yet the Babylonians computed three meaningful rotation rates for the moon's orbit moving on a three-dimensional torus.
\allblack

Our approach is to project the torus into the plane or onto a circle 
$\phi:\torus\to S^1$
and define the rotation rate as the rate of rotation of the image.

\allblack
\subsection{ Rotation rate of a torus defined relative to 
$\BF{\phi:\torus\to S^1}$} \label{ss:phi}
For completeness, we reproduce some material from our paper \cite{DCDS} since without this material, it is not clear why we always define a rotation rate relative to a map $\phi:\torus\to S^1$. Such maps  have a nice representation.
Let $a = (a_1,\cdots,a_d)$ where $a_1,\cdots,a_d$ are integers and let 
$\theta = (\theta_1,\cdots,\theta_d)\in\torus$. 
The simplest $\phi$ has the form
$\phi(\theta) := a_1\theta_1+\cdots+a_d\theta_d \bmod1,$ which can also be written $e^{i2\pi(a_1\theta_1+\cdots+a_d\theta_d)}\in S^1$ if $S^1$ is viewed as the unit circle in $\mathbb{C}$. Then $\phi$ is a continuous map of the torus to a circle.
For any initial point $\theta_0 \in \torus$, we have $\theta_n =  \theta_0 + n(a_1\rho_1+\cdots+a_d\rho_d )\bmod1$ and in this very simple case 
$ \theta_{n+1}-\theta_n = a\cdot\rho \bmod1:= a_1\rho_1+\cdots+a_d\rho_d \bmod1$ is constant and in this very special case we obtain a constant rotation rate for $\phi(\theta)$, namely

\begin{equation}\label{eqn:rho_phi}
\rho_\phi \bmod1=a\cdot\rho \bmod1.
\end{equation}
See Eq. \ref{eqn:a_dot_rho}. 
For $d=1$, Eq. \ref{eqn:rho_phi} says $\rho_{\phi}=a_1 \rho\bmod1$ where $a_1$ is an integer. 
The integer $a_1$ depends on the choice of $\phi$, so even when $|a_1| = 1$ we can  
get $\rho$ for one choice and $-\rho$ for another choice.

\begin{equation}\label{eqn:rr}
\rho_\phi :=  a\cdot \rho\bmod1.
\end {equation}
We note that for every map $\phi$ of a torus to a circle, there are integers $a_j$ and a periodic function $g:\torus\to\R$ such that
\begin{equation}\label{eqn:g}
\phi(\theta) = g(\theta) +a\cdot \theta \bmod1.
\end{equation}
Directly computing a rotation rate for this map can be difficult. 
In fact, after we define the rotation rate below, it will turn out that Eq. \ref{eqn:rr} will still be true, independent of $g$, but this formula will not be very helpful in determining $\rho_\phi$ from the image of a trajectory, $\phi(\theta_n)$ since $a$ is unknown.

\textbf{Defining a rotation rate for  
 $\phi:\torus\to S^1$.}
Rotation rates are key characteristics of any quasiperiodic
trajectory. One heuristic way of thinking of the rotation rate of $\phi$ for a trajectory $(\theta_n)$ on $\torus$ is to say that it is the average of the angle differences $\theta_{n+1}-\theta_n$, but one cannot average angles unless these differences are roughly constant. Then the angle differences can be continuously ``lifted'' to the real line, and the average is a real number.
 We now make this more precise. 
\allblack
%\allred

%\textbf{Defining $\Delta$ and its lift $\hat\Delta$ for  $\phi:\torus\to S^1$.}
%
Suppose 
there exists a continuous 
map $\phi:\torus\to S^1$ from the dynamical
system to a circle, 
%and we have available 
but we only know 
the image
$\phi_n:=\phi(n\rho)$  
sequence of a trajectory $F(\theta_n)=
\theta_{n+1}=\theta_{n}+\rho\bmod1$ on a torus.
%\newpage
Define 
\begin{align}
\Delta(\theta) :&= \phi(\theta+\rho)-\phi(\theta)\bmod1\nonumber\\
 &= g(\theta+\rho)+ a\cdot (\theta+ \rho) -[ g(\theta)+ a\cdot (\theta)]\bmod1\mbox{ (from Eq. \ref{eqn:g})}\nonumber\\
 &= a\cdot \rho + g(\theta+\rho) -g(\theta) \bmod1.\label{Eq5}
\end{align}
We say $\hat\Delta$ is a {\bf lift} of $\Delta:\torus\to S^1$ if 
(i)$\hat\Delta:\torus\to\mathbb{R}$, (ii)
$\hat\Delta$ is continuous; and (iii) $\hat\Delta(\theta)\bmod1=\Delta(\theta)$. 
Motivated by Eq. \ref{Eq5}, we define
 \begin{equation}\label{hatdelta}
   \hat\Delta(\theta) :=a\cdot \rho + g(\theta+\rho) -g(\theta)\in\R.
\end{equation}
Then (i),(ii), and (iii) are satisfied so $\hat\Delta$ is a lift of 
$\Delta$. Define $\hat\Delta_n= \hat\Delta(\theta_n)$. 
We can formally define the {\bf rotation rate for} $\mathbf{\phi}$ as
\begin{equation}\label{eqn:delta}
\rho_\phi := \left(\displaystyle\lim_{N\to\infty}\frac{\sum_{n=0}^{N-1}\hat\Delta_n}{N}\right)\bmod1.
\end{equation}
\begin{prp}
Assume $\theta_n$ is quasiperiodic.
then the limit $\rho_\phi$ in Eq.~\ref{eqn:delta} exists,
and 
\begin{equation}\label{eqn:a_dot_rho}
\rho_\phi = a\cdot\rho\bmod1.
\end{equation}
\end{prp}

\begin{proofn}
The existence of the limit is guaranteed by the Birkhoff Ergodic Theorem, which says that the limit in Eq. \ref{eqn:delta} is
\begin{align}
\int_\torus \hat\Delta(\theta)d\theta  &=
\int_\torus\big(
a\cdot \rho + g(\theta+\rho) -g(\theta)\big)d\theta\nonumber \\
&=a\cdot \rho + \int_\torus g(\theta+\rho)d\theta -\int_\torus g(\theta)d\theta\label{eq7}\\
&=a\cdot \rho,
\end{align}
since $\int_\torus d\theta = 1$ and the two integrals of $g$ in Eq. \ref{eq7} are equal.
Hence $\rho_\phi =a\cdot \rho\bmod1.$\qed
\end{proofn}

Different choices of the lift $\hat\Delta$ can change $\rho_{\phi}$ by an integer, so $\rho_{\phi} \bmod1$ is independent of the choice of lift $\hat\Delta$.
The {\bf rotation rate} is this $\rho_\phi \bmod1$. 
The limit in Eq. \ref{eqn:delta} exists and is the same for all initial $\theta_0$. 

\allblack

\section{Numerical results for a one-dimensional example 1 in Eq. \ref{eq:1d}}\label{s:siegel-disk}
\allblack
\subsection{Computing a conjugacy}
Fourier series has long played an important role in investigations of quasiperiodicity. See for example
\cite{gomez_2010a,gomez_2010b}.%,sanchez_2010}.
\allblack

Of course numerically we only determine the Fourier coefficients $b_k$ in 
Eq.~\ref{eq:H} 
for $|k| \le K$, stopping when it appears that  
$|b_k| < 10^{-30}$ for larger $k$. 
We determine $b_k$ with an error of less than $10^{-30}$. 
%\allred
But what about the accuracy of the trajectory, which requires iteration of the map? If we iterate the linear map and take the conjugacy of each point, we obtain a pseudo-trajectory that we denote by $\hat z_n$. We would like it to be close to $z_n$,  the trajectory  obtained by numerically iterating the map. The difference will depend on the trajectory, so here we have chosen our most nonlinear trajectory  called ``magenta'' in Fig.~\ref{fig:quasi-siegel}~(Right).  

We find below that the trajectory produced from the Fourier series differs from by at most $3\cdot10^{-26}$.

Write $\hat b_k$ for $|k| \le K$ and $\hat\rho$ for our computed values of $b_k$ and $\rho$.
We can computationally reconstruct $z_n$ by  
$$
\hat z_{ n}= h(n\hat\rho)
=\sum_{k=0}^{K} \hat b_k \sigma_k(n\hat\rho).
$$

How close are $\hat z_{ n}$ to $z_n:=F^n(z_0)$?
\allblack
For Eq. \ref{eq:1d} with initial point 
%\allred
$z_0 = 0.37 \in \mathbb{C}$ we use  iterates $z_n$ for $n=0,\cdots,N$, where $N=4\cdot10^6$
%\allblack
 to compute $K:=3400$ Fourier coefficients. We find
$$
|z_n - \hat z_n| < 3\cdot10^{-26}\mbox{ for all }n\le N.
$$
%\allblack

%
One can compute the power series of the Siegel conjugacy function iteratively from a knowledge of $F$; see Problem 11-c in \cite{Milnor},
%Milnor's complex book,
 a method discussed by Cremer~\cite{Cremer}; 
 in iterative schemes, the error can increase and it is necessary to control the cumulative error. 
 There is also a way of computing conjugacy using Newton method to find Fourier coefficients. 
See for example~Jorba~\cite{NumericQuasi5}.
 
 We prefer our direct computation of the Fourier coefficients and we do not need access to $F$.
%\allred
%{\bf Our main example in $\mathbf{\mathbb{C}}$, continued.}
%\allblack
Fig.~\ref{fig:quasi-siegel} shows several closed curves, each of which is invariant for Eq.~\ref{eq:1d} and are quasiperiodic for 
$$
z_{n+1}= F(z_n).
$$
First we focus on four of the curves except for the outermost curve.
Although the inner curve shown looks like a smooth circle, the outer curve plotted looks rough-edged. 
All are analytic. 
%\allred
Each curve is conjugate to 
%has a highly nonlinear choice of variables 
$\theta$ so that 
\allblack
$$
\theta_{n+1}=\theta_n+ \rho.
$$

Rotation rates for the inner four orbits are calculated by employing 
%a weighted Birkhoff average 
$\Q$ 
%method 
explained in Sec.~\ref{ss:WB}.
Since we know the actual value of $\rho$, we can determine the error in our calculations.
For each curve, $\rho$ is computed to a precision of about $10^{-33}$ using
$N$ iterates of Eq.~(\ref{eq:1d}) 
%\allred
%where $N \approx 3000, 10^5,$ and $10^6$ for the four curves used
 in Figs. \ref{fig:convergence-siegel}.
 %and \ref{fig:fourier-siegel}.
%\allblack
Fig.~\ref{fig:convergence-siegel} shows that the computed rotation rate $\rho_N$ converges to the rigorous rotation rate $\rho$ as the number of iterations $N$ increases, 
and the errors in the rotation rate $err(N)=|\rho_N -\rho|$.

%\allred
%for the four curves reach about $10^{-32}$ when $N=8000, 15000, 20000, 300000$, respectively.
%\allblack

\allblue
In Figs.~\ref{fig:fourier-siegel}, \ref{fig:fourier-siegel-converted}  the four curves investigated have 
initial values $z_0= 0.1, 0.2, 0.3$, and $0.37$.
Fig.~\ref{fig:fourier-siegel-converted} (Left) and Fig.~\ref{fig:fourier-siegel}  show the magnitude of the Fourier coefficients of $h$, the change of coordinates 
$\phi$ between the quasiperiodic circle and the pure rotation with the rotation rate $\rho$.
Although the rate of decrease in the magnitudes is slow for the outer curve, 
all of the four curves show the exponential decay at least within the quadruple precision accuracy.   
Only from the knowledge of the trajectory  we can obtain a strong implication that 
%These strongly suggest that
 such a curve is computationally analytic.
The slopes are the values of $R_0$. Their values reported in those log-log plots are
$R_0=0.285,0.570,0.836,0.984$, respectively.

In Fig.~\ref{fig:convergence-siegel} we can see that to get an accuracy of $10^{-30}$, we require increasing values of $N$, and as $1-R_0\to 0,$
the $N$ seems to diverge. In the graph we can see that the required $N$ is approximately
$N=2000, 3000, 10^4, 10^5$, respectively for the four cases shown. 
These numbers are roughly proportional to $(1-R_0)^{-1}$.
% where the $R_0$ values are shown in the previous two figures, Figs.~\ref{fig:fourier-siegel-converted},\ref{fig:fourier-siegel}. Specifically $(1-R_0)^{-1}= 1.4,2.3,6,60$
 Specifically
$1666\cdot(1-R_0)^{-1}\approx
2300,3900,10^4,10^5$.
\allblack

%\allblue
\subsection{Estimation of the $\mathbf{L^2}$ length of the quasiperiodic curve
$\mathbf{h(C_r)}$}
\begin{figure}
  \centering
    \includegraphics[width = .4\textwidth,height=.4\textwidth]{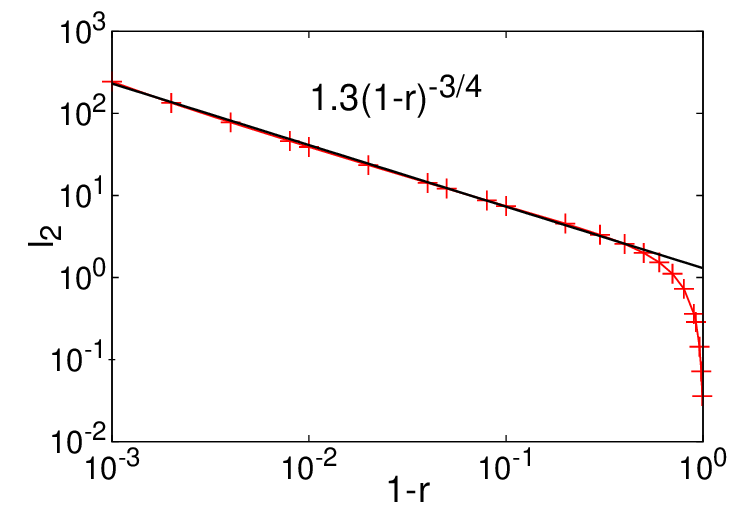}
  \includegraphics[width = .4\textwidth,height=.4\textwidth]{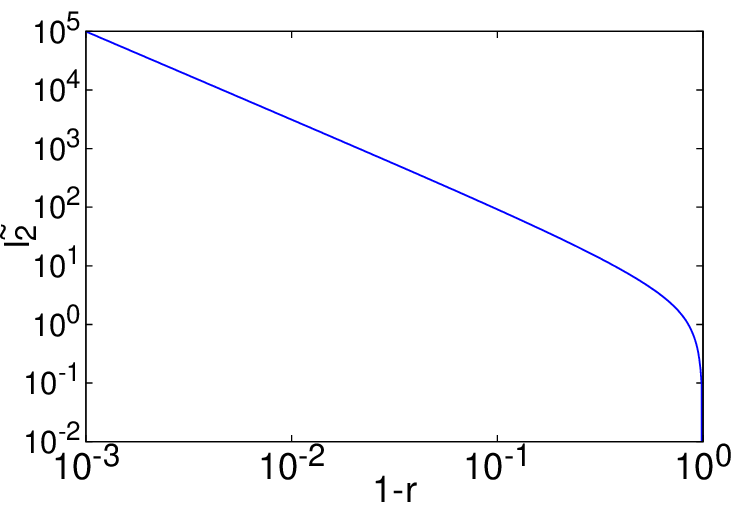}
  \caption{
%\allred
  \textbf{Estimates of the total $L^2$  length of $h(C_r)$ in the Siegel disk.} 
    Left: $L_2$ length $l_2(\gamma)$ is computed using the direct calculation of Eq.~\ref{eq:L2} with the estimated value of $a_k$ created in  Fig.~\ref{fig:fourier-siegel-converted} for some values of $r$, which approaches to $1.3(1-r)^{-3/4}$, as $r\to 1$.
    Right: The upper bound $\tilde{l_2}(\gamma)$ in Ineq.~\ref{ineq:length} is shown for $c=1.$
    %using Ineq.~\ref{ineq:length}.
    }
  \label{fig:length-siegel}
\end{figure}

% in the Siegel disk.} 
As $r\nearrow 1$, the curve $h(C_r)$ becomes even more irregular, while remaining analytic.
Let $\gamma:[0,1)\to\mathbb{C}$ be $C^1$. Then we define the $L^2$ length of $\gamma$ to be
$$
	l_2(\gamma) := \left[ \int_0^1~\left|\frac{d}{d\theta}\gamma(\theta)\right|^2d\theta\right]^{1/2}.
$$
\begin{prp}
Let $\D$ be the unit disk, and  $h:\D\to\mathbb{C}$ be a linearization conjugacy in Eq.~\ref{eq:conjugacy}.
Assume the power series for $h(z)=\displaystyle\sum^{\infty}_{k=0}  a_kz^k$ has the property that
$c := sup~|a_k| < \infty$. Then
there exists $c>0$ such that $|a_k|\le c$ for all $k\ge 0$.
Let $0\le r <1$.  
Define $\gamma(\theta):=h(re^{i 2\pi\theta}) $ for $\theta\in[0,1)$.
Then 
\begin{align}
l_2(\gamma)
&= 2\pi\left( \sum^{\infty}_{k=0} \left|ka_kr^k\right|^2\right)^{1/2}\label{eq:L2}\\
&\le2\pi c\left[\frac{r^2(1+r^2)}{(1-r^2)^3}\right]^{1/2}.\label{ineq:L2}
\end{align}
If $|a_k|= c$ for all $k\ge 0$, then the inequality in Ineq.~\ref{ineq:L2} becomes an equality.
\end{prp}
%
%$$z=re^{i 2\pi\theta}$$
%
\allblack
%{\bf Proof.} 
\begin{proofn}
The assumptions imply
$$
\gamma(\theta)= \sum^{\infty}_{k=0} a_kr^k\sigma_k(\theta)\
%mbox{ and}
$$
and
\begin{equation}
\frac{d}{d\theta}h(re^{i2\pi\theta})
%\frac{dH}{d\theta}(\theta) 
= 2\pi i \sum^{\infty}_{k=0}  ka_kr^k\sigma_k(\theta),
\end{equation}
where  $\sigma_k(\theta)=e^{2\pi i k \theta}$.
%\allblack
The $\sigma_k$ are orthogonal in $L^2[0,1)$ and 
$|\sigma_k(\theta)|=1$ and 
$|a_k| \le c$, so
\begin{align}%\label{eq:l2square}
l_2(H)^2
&=\int_0^1~\left|\frac{d}{d\theta}h(re^{i2\pi\theta})\right|^2d\theta 
= (2\pi)^2  \sum^{\infty}_{k=0}  \int_0^1~\left|ka_kr^k\sigma_k(\theta)\right|^2 d\theta
= (2\pi)^2  \sum^{\infty}_{k=0} \left|ka_kr^k\right|^2\label{eq:length}\\
&\le (2\pi c)^2  \sum^{\infty}_{k=0}  ~|kr^k|^2=:\tilde{l_2}(H)^2.\label{ineq:length}
\end{align}
We can rewrite Eq.~\ref{ineq:length} by setting $\psi=r^2 (<1)$ and obtain $(2\pi c)^2  \displaystyle\sum^{\infty}_{k=0}  ~k^2\psi^k$. 
Using a well known formula, 
% (see for example, 
%\begin{verbatim}
%Wikipedia (https://en.wikipedia.org/wiki/List_of_mathematical_series)), 
%\end{verbatim}
\begin{equation}
\displaystyle\sum^{K}_{k=1}k^2 \psi^k
=\psi\frac{1+\psi
\red{-(K+1)^2 \psi^K+(2K^2+2K-1)\psi^{K+1}-K^2 \psi^{K+2}}}{(1-\psi)^3}\label{eq:l2square}
%\to \frac{\psi(1+\psi)}{(1-\psi)^3}~\text{as}~K\to\infty,
\end{equation}
for $|\psi|<1$, in the limit as $K\to\infty$ the colored terms go to $0$, so
we can compute the final sum in Eq.~\ref{eq:l2square}
\begin{equation}
	\displaystyle\sum^\infty_{k=0} k^2 \psi^k
	=\frac{\psi(1+\psi)}{(1-\psi)^3}
=\frac{r^2(1+r^2)}{(1-r^2)^3}.
\end{equation}
Then we obtain Eq.~\ref{eq:L2} and Ineq.~\ref{ineq:L2}.\qed
\end{proofn}
%Then $l_2(H)$ is obtained as a function of $r$.
%\qed\\
We have included the terms in red for the readers convenience, because they can be useful in making numerical computations.\\
In Fig.~\ref{fig:fourier-siegel-converted}, right panel, we can see that for our case the $|a_k|$ are far from constant so that the formula in Eq.~\ref{eq:L2} can be more useful than that in Ineq.~\ref{ineq:L2}.
Aside from the initial $k < 1000$, the values of $|a_k|$ are bounded above by $10^{-4}$. 
Then the estimate in Ineq.~\ref{ineq:L2} is not sharp. 
%Then the estimate in Ineq.~\ref{ineq:length} is not sharp. 
Therefore, in Fig.~\ref{fig:length-siegel}~(Left) we compute the values of $l_2(\gamma)$ directly using the observed values of $a_k$ for various values of $0<r<1$, in comparison with $\tilde{l_2}(\gamma)$.

%
%\ifincludeXX
%\fi
%\newpage

\section{Numerical results for a two-dimensional example 2 in Eq. \ref{henon}: H\'enon map}\label{s:siegel-ball}
%\section{Quasiperiodic orbit within a Siegel ball in the complex H\'enon map: Example 2}\label{s:siegel-ball}
%{\bf Complex H\'enon map: Example 2.}
%
Siegel ball is the analogue of the Siegel disk for a higher dimensional complex dynamical system.
In this section, we compute rotation vectors of a quasiperiodic curve within a Siegel ball for the complex H\'enon map defined in  Eq.~\ref{henon}. See \cite{Ushiki0}.
%\allred
%That is, SIEGEL BALL 
%\allblack

{\bf Choosing $\alpha$ and $\beta$ for the H\'enon map Eq.~\ref{henon}.}
\allblack
We specify $\alpha$ and $\beta$ in terms of parameters $\theta, \phi \in \R$:
\begin{align}
  \alpha&=2\cos \theta \cos \phi -\cos^2 \phi,\label{alpha}\\
  \beta&=e^{i\theta}. \label{beta}
%  \beta&=\cos \theta +i \sin \theta. 
\end{align}
Fixed points of $F_2:=F_{\alpha,\beta}$ are $y_{*}=\cos \theta \pm (\cos \theta -\cos \phi).$
We focus on $y_{*}=\cos \phi$ and set $\mu=\cos \phi + i \sin \phi$.
The Jacobian matrix at the fixed point is 
\begin{equation}
DF_{\alpha,\beta}=\beta 
    \begin{pmatrix}
      0 & 1/\beta\\
      -\beta & 2\cos \phi
    \end{pmatrix}.
\end{equation}
Note that it has determinant $+1$.
Its eigenvalues  
are $\beta \mu$ and 
$\beta \overline{\mu}$;
in other words $(e^{i[\theta\pm\phi]})$.
The corresponding eigenvectors are 
$$
(e^{\mp i\phi},e^{i\theta})\in\mathbb{C}^2.
$$
\allblack
The linearized system is quasiperiodic (for appropriate parameter choices) and we have to say the nonlinear system has the same rotation vectors. 
Therefore the rotation vector obtained through a projection will be expressed by some linear combination of 
$\rho_1=(\theta + \phi)/2\pi$ ~and~ $\rho_2=(\theta-\phi)/2\pi$.
%implying the rotation rates are 
%\begin{equation}
% (\theta + \phi)/2\pi \hspace{1zw} \text{and} \hspace{1zw} (\theta-\phi)/2\pi.
%\end{equation}
%At least the former eigenvalue is confirmed to be related to the rotation rate of an orbit $((\theta + \phi)/2\pi)$.

\subsection{Example 2A}
%%\ifincludeXX
\begin{figure}[t]
\centering{
  \includegraphics[width= .46\textwidth,height=.35\textwidth,angle=0]{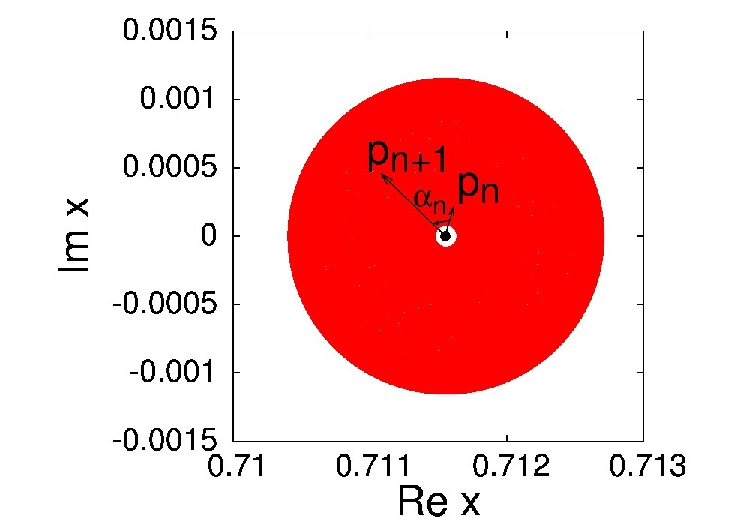}
  \includegraphics[width= .46\textwidth,height=.35\textwidth,angle=0]{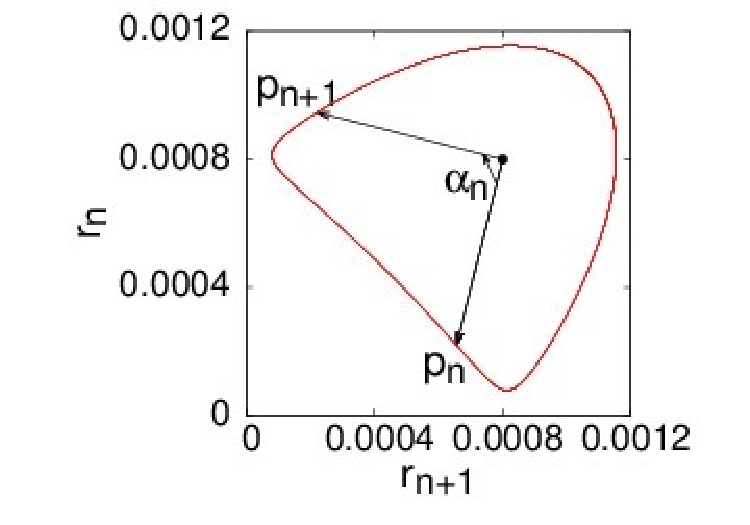}
}
\caption{{\bf Projections of a trajectory of Example 2A.}
In order to compute a rotation rate we project the torus into the plane. 
If the projection of the torus to the plane has a hole in it as seen in both panels, we can compute the rotation number relative to a reference point in that hole (and as we mention elsewhere) the resulting rotation number is the same for all choices of reference point. In both panels, $\alpha_n$ is the angle between consecutive points. But $\alpha_n\in S^1$, or in 
$[0,1)\bmod1$.
Division is not uniquely defined in $S^1$. 
Therefore we cannot average the values $\alpha_n$. As discussed in the Introduction, it is necessary to create a continuous lift of $\alpha$ to $\R^1$ to obtain numbers that can be averaged.
The setting for both panels: $\rho_1=(\sqrt{5}-1)/2=0.61803398875\cdots,\rho_2=\sqrt{3}/2=0.86602540378\cdots,\theta=(\rho_1 +\rho_2)\pi, \phi=(\rho_1 - \rho_2)\pi$. 
Left: Projection of an orbit onto $(\text{Re}~x, \text{Im}~x)$ plane. 
Right: Projection of an orbit onto ($r_{n+1},r_{n}$), where $r=\sqrt{(\text{Re}~x-u)^2 +(\text{Im}~x -v)^2}$ where 
$(u,v)$ is the reference point $(\cos(\phi),0)=(0.71155\cdots,0)$ used in the left panel.
}
\label{ex1-proj}
\end{figure}
%%\fi

\begin{figure}[hbtp]
\centering{
\includegraphics[width= .42\textwidth,height= .32\textwidth,angle=0]{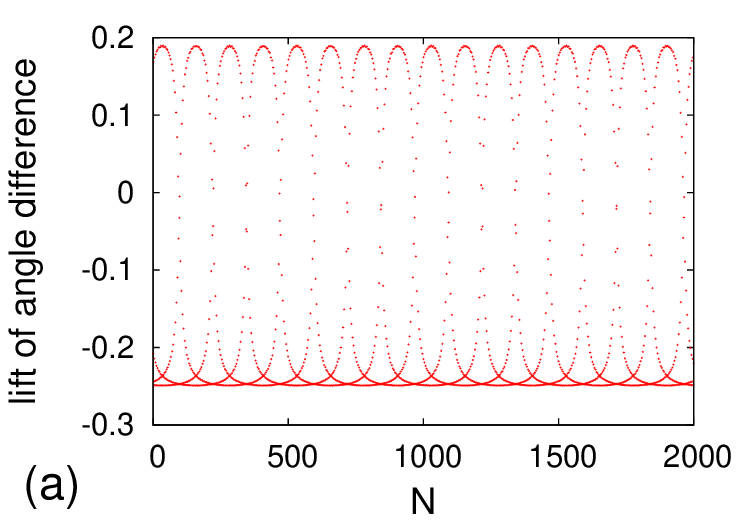}
\includegraphics[width= .42\textwidth,height= .32\textwidth,angle=0]{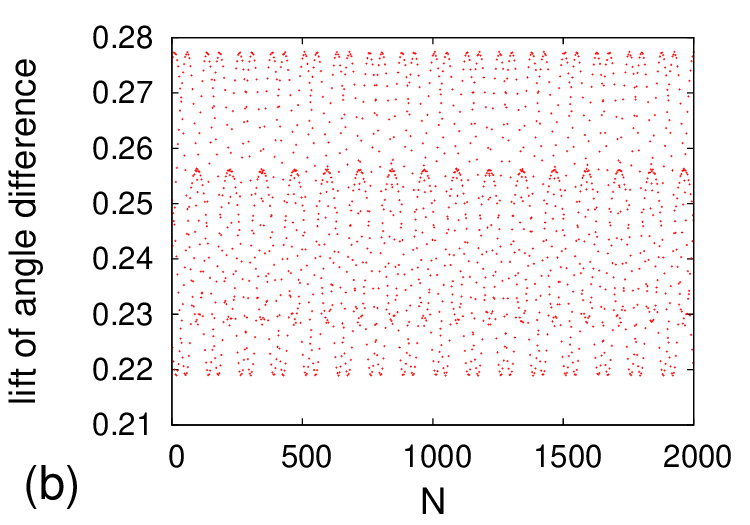}
\includegraphics[width= .42\textwidth,height= .32\textwidth,angle=0]{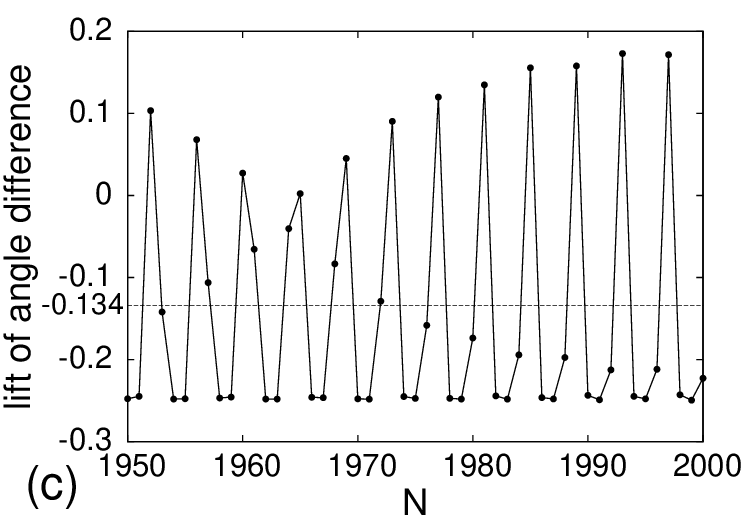}
\includegraphics[width= .42\textwidth,height= .32\textwidth,angle=0]{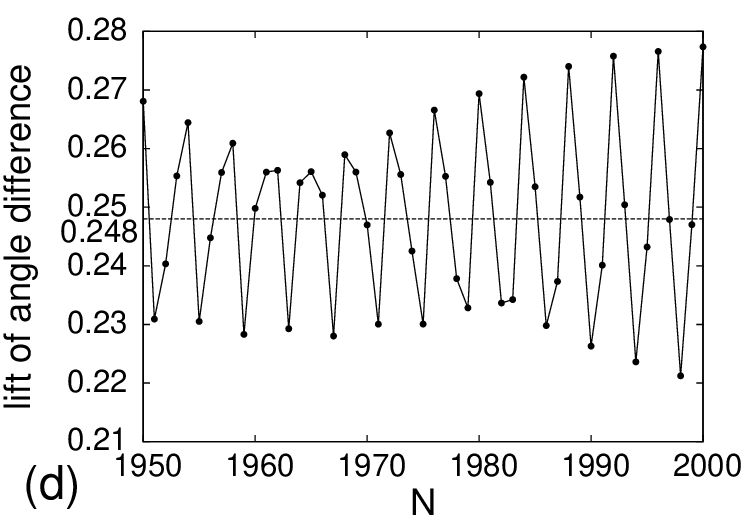}
\includegraphics[width= .42\textwidth,height= .32\textwidth,angle=0]{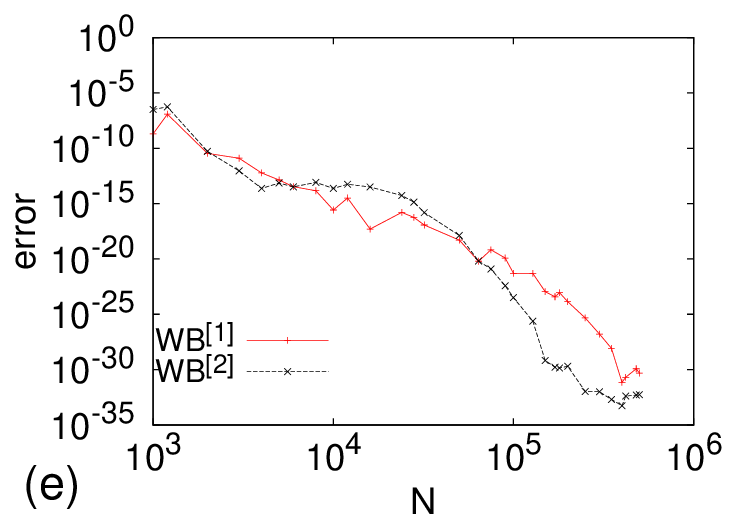}
\includegraphics[width= .42\textwidth,height= .32\textwidth,angle=0]{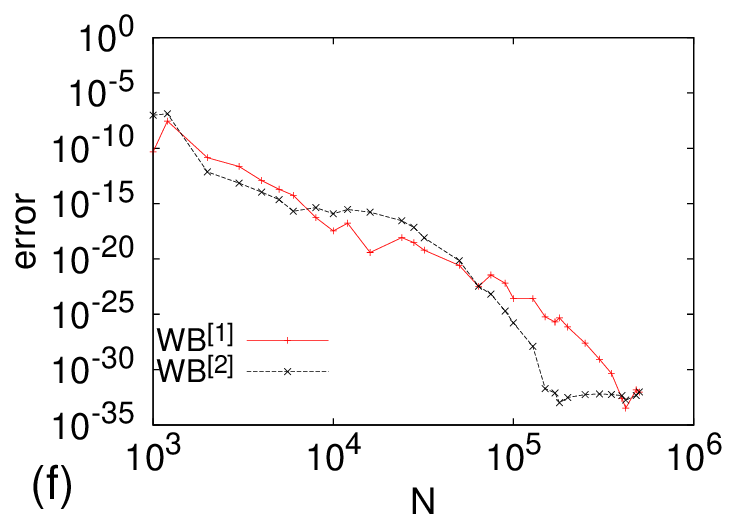}
}
\caption{{\bf Computation of rotation rates of Example 2A.} (a) Lift of the angle difference about Fig.~\ref{ex1-proj}~(Left) viewed from $(u, v)=(\cos(\phi),0)=(0.71155\cdots,0)$. (b) Lift of the angle difference about Fig.~\ref{ex1-proj}~(Right) viewed from $(r_{n+1},r_n)=(0.0008,0.0008)$.
(c) and (d) Detail versions of  (a) and (b) with consecutive dots connected.
(e) The Birkhoff average ($\Q^{[1]}$ and $\Q^{[2]}$) of the lift of the angle difference shown in (a) has converged to $-0.13397~45962~15561~35323~62768~29247
%~06685$ 
\sim \rho_2-1$, and the error was $\sim10^{-30}$. (f) That in (b) has converged to $0.24799~14150~34543~79855~91363~36387
%~30229$ 
\sim  (1-\rho_1)-(1-\rho_2)$, and the error was $\sim10^{-32}$.
}
\label{ex1-conv}
\end{figure}

We investigate a quasiperiodic orbit within a Siegel ball for a set of parameter values $\theta=(\rho_1 +\rho_2)\pi$ and 
$\phi=(\rho_1 - \rho_2)\pi$, where $\rho_1=(\sqrt{5}-1)/2$ and $\rho_2=\sqrt{3}/2$. Fig. \ref{ex1-proj}~(Left) shows an orbit 
projected onto $(\text{Re}~x,\text{Im}~x)$ plane, and Fig.~\ref{ex1-proj}~(Right) shows the same orbit projected onto  
($r_{n+1},r_{n}$) plane, where $r=\sqrt{(\text{Re}~x-u)^2 +(\text{Im}~x -v)^2}$ and $(u,v)=(\cos(\phi),0)=(0.71155\cdots,0)$.
The lift of the angle difference about Fig.~\ref{ex1-proj}~(Left) viewed from $(u, v)=(\cos(\phi),0)$ is shown in Fig.~\ref{ex1-conv}~(a) and its detail in (c), and that about Fig.~\ref{ex1-proj}~(Right) is shown in Fig.~\ref{ex1-conv}~(b) and its detail in (d).
Fig.~\ref{ex1-conv}~(e) and (f) show the convergence of the Birkhoff average of the lift of the angle difference (a) and (b), respectively.
The Birkhoff average is calculated using the weighted Birkhoff average with two different weight functions ($\Q^{[1]}$ and $\Q^{[2]}$) described in Sec.~\ref{ss:WB}.
% in two ways, the weighted Birkhoff average ($\Q^{[1]}$), the weighted Birkhoff average with a different weight function ($\Q^{[2]}$).
%, and the original Birkhoff average ($\B$). 
The faster convergence in the quadruple precision is realized when we use $\Q^{[2]}$; see \cite{DCDS}.
The computed rotation rate in Fig.~\ref{ex1-conv}~(c) is $-0.13397~45962~15561~35323~62768~29247$
%~06685$ 
%
$\sim\rho_2-1$,
%=$\sqrt{3}/2-1$, 
and that obtained in Fig.~\ref{ex1-conv}~(d) is $0.24799~14150~34543~79855~91363~36387$
%~30229 
$\sim\rho_2-\rho_1$.
%=$\sqrt{3}/2-(\sqrt{5}-1)/2$.
Hence, both of these rotation rates are integer-sum combinations of $\rho_1$ and $\rho_2$ $\bmod1$.
\allblack
See the explanation in Sec.~\ref{ss:phi}.
%%\ifincludeXX
%%\fi

\newpage
\subsection{Example 2B}\label{ss:2B}
\begin{figure}[h]
\centering{
  \includegraphics[width=.35\textwidth,height=.32\textwidth,angle=0]{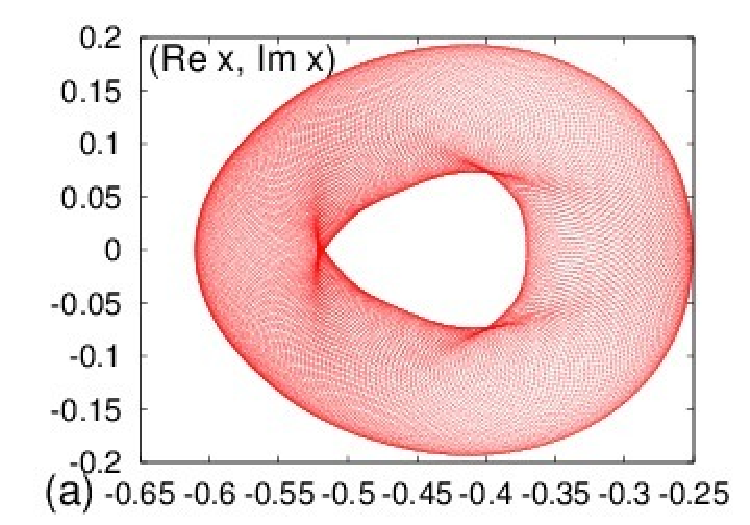}
  \includegraphics[width=.35\textwidth,height=.32\textwidth,angle=0]{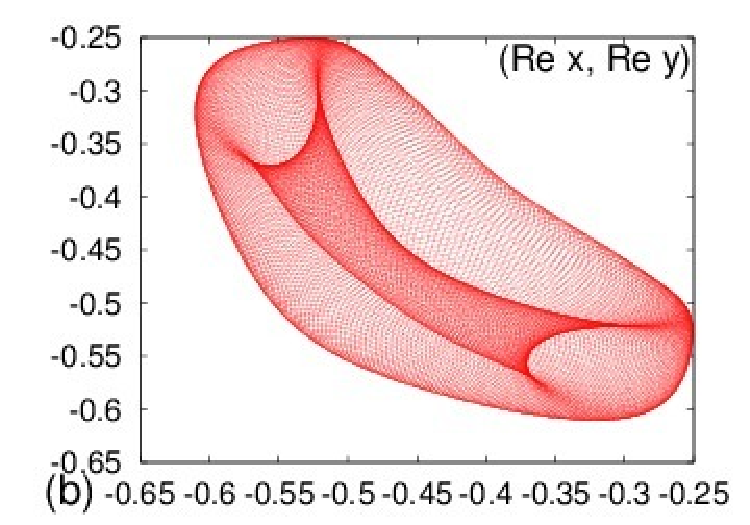}
  \includegraphics[width=.35\textwidth,height=.32\textwidth,angle=0]{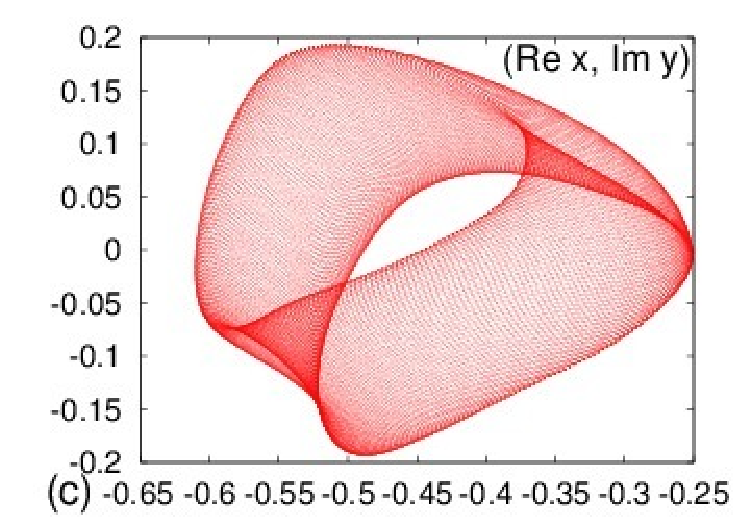}
  \includegraphics[width=.35\textwidth,height=.32\textwidth,angle=0]{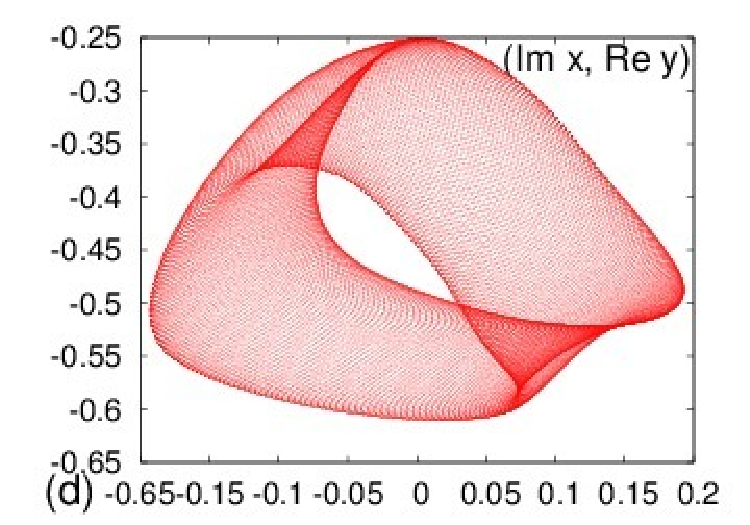}
  \includegraphics[width=.35\textwidth,height=.32\textwidth,angle=0]{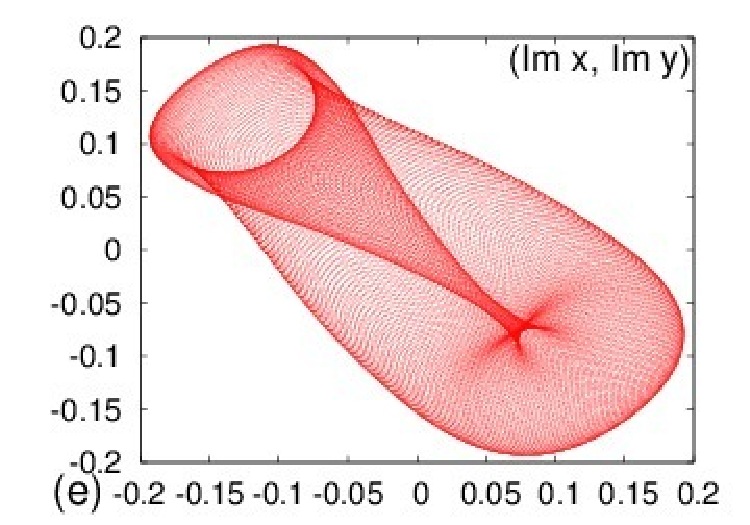}
  \includegraphics[width=.35\textwidth,height=.32\textwidth,angle=0]{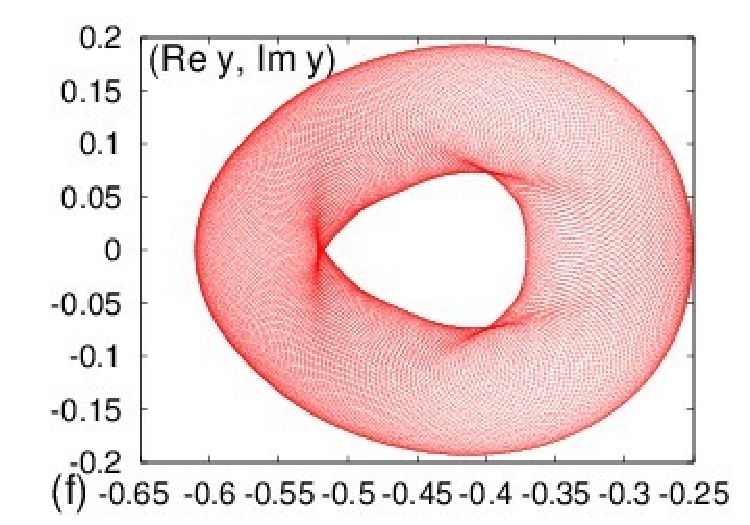} 
}
\caption{{\bf Different projections of a trajectory of Example 2B.} The initial point is $x_0=-0.500+0.126i, y_0=-0.387-0.163i$, 
where the parameters are set as $\theta=0.664, \phi=2.032$ in Eqs.~\ref{alpha} and \ref{beta}.
In (b) and (e) there is no hole in the projection so a rotation rate cannot be computed.
In (a),(c),(d),(f) there is an interior white region about which a rotation rate can be computed. See Sec.~\ref{ss:2B}.
}
\label{simple}
\end{figure}

\begin{figure}[h]
\centering{
    \includegraphics[width=7.5cm,height=5.5cm,angle=0]{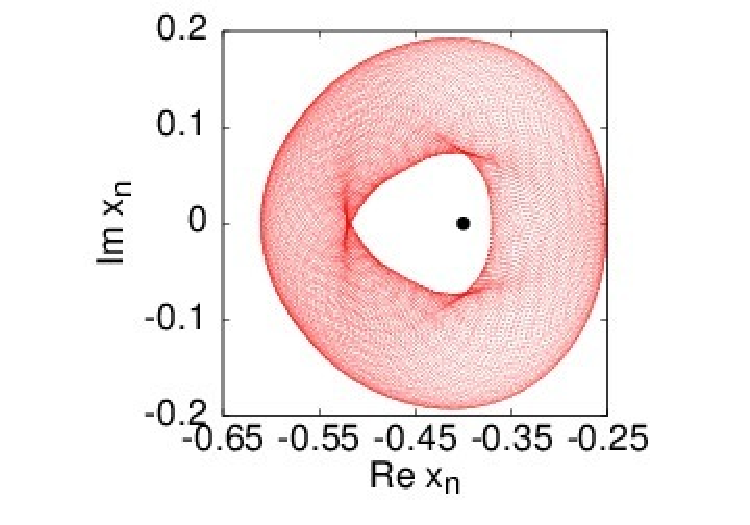}
       \includegraphics[width=7.5cm,height=5.5cm,angle=0]{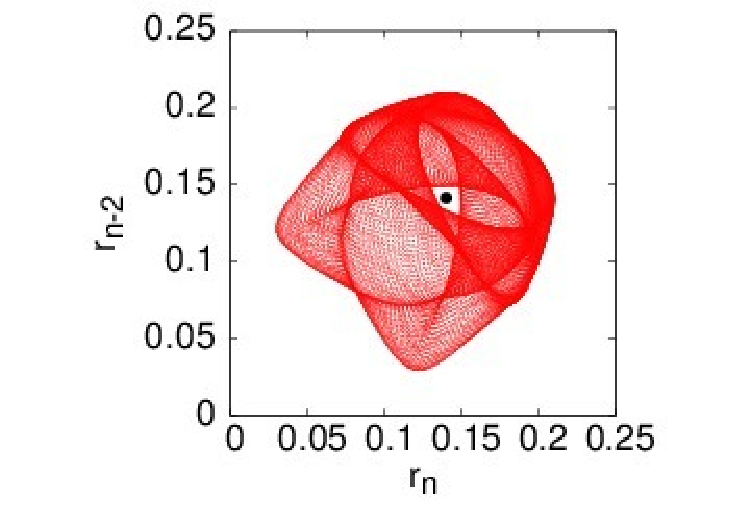}
}
\caption{{\bf Projections of a trajectory onto two different planes of Example 2B.} Left: $(\text{Re}~x,\text{Im}~x)$ plane.  
%started from a point $(x_0,y_0)=(-0.500+0.126i,-0.387-0.163i)$ when parameters are set as $\theta=0.664, \phi=2.032$, 
Right: $(r_n,r_{n-2})$ plane, where $r_n=\sqrt{(\text{Re}~x_n +0.4)^2 + {\text{Im}~x_n}^2}$.
In each panel the black point in the central region 
($(\text{Re}~x,\text{Im}~x)=(-0.4,0)$ on the Left and $(r_n,r_{n-2})=(0.14,0.145)$ on the Right)
indicates the reference point used to measure the rotation rate.
Moving the black dot within the central white region makes no difference in the resulting rotation rate. See Fig.~\ref{ex2-conv} for 
the computation of rotation rates.
%numerical results.
}
\label{ex2-proj}
\end{figure}

\begin{figure}[h]
\centering{
\includegraphics[width=7.5cm,height=5.75cm,angle=0]{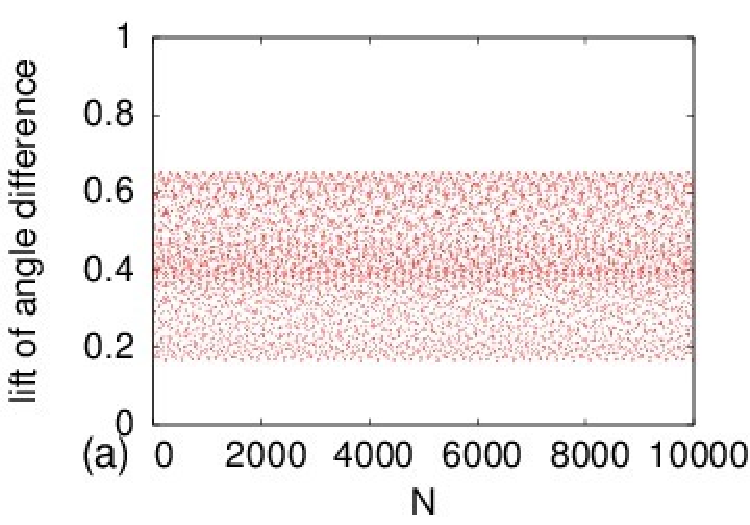}
\includegraphics[width=7.5cm,height=5.75cm,angle=0]{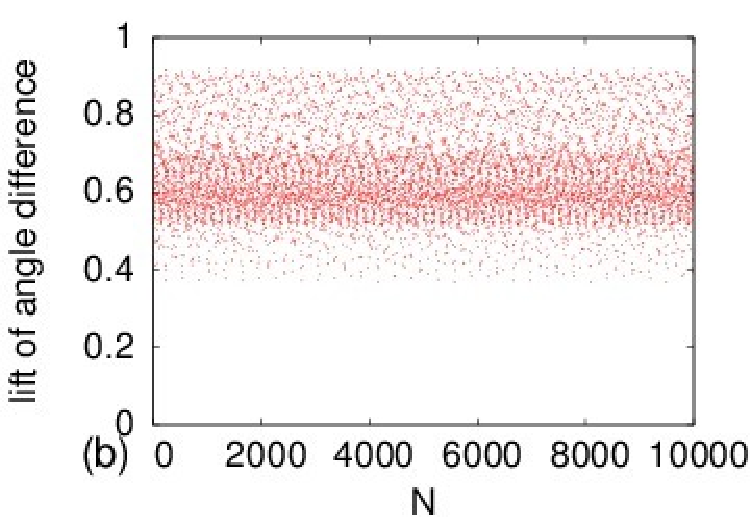}
\includegraphics[width=7.5cm,height=5.75cm,angle=0]{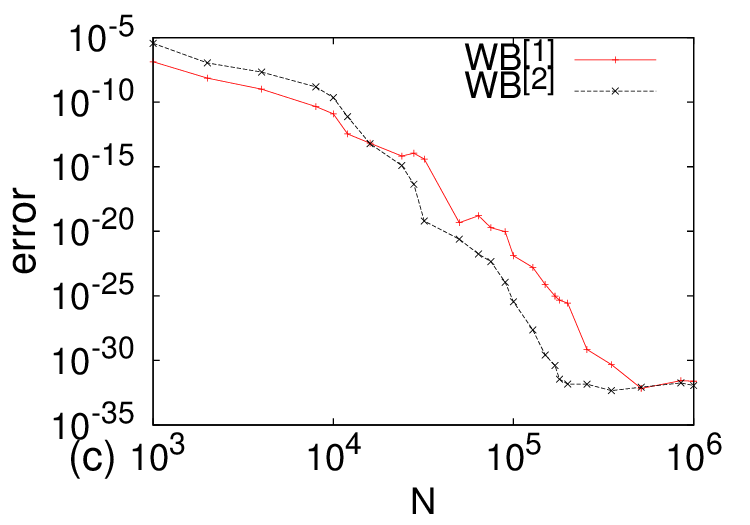}
\includegraphics[width=7.5cm,height=5.75cm,angle=0]{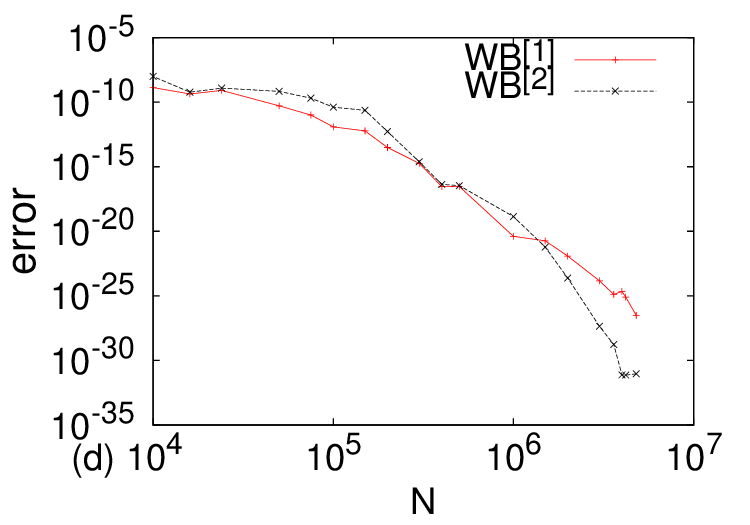}
}
\caption{{\bf Computation of rotation rates of Example 2B.} (a) The lift of angle difference for Fig.~\ref{ex2-proj}~(Left).
(b) The lift of the angle difference for Fig.~\ref{ex2-proj}~(Right).
% and (b) $(r_n,r_{n-5})=(0.135,0.135)$
%, where $r_n=\sqrt{(\text{Re}x_n +0.4)^2 + {\text{Im}x_n}^2}$.
(c) The Birkhoff average (using $\Q^{[1]}$ and $\Q^{[2]}$) of the lift of the angle difference in (a) has converged to  $\rho_1=(\theta+\phi)/2 \pi\sim 0.42908~17265~75749~82523~29106~26052$.  The error is $\sim 10^{-31}.$
(d) 
%The Birkhoff average ($\Q^{[1]}$ and $\Q^{[2]}$) of the lift of angle difference 
For (b)  $\Q^{[2]}$  has converged to $\rho_1-\rho_2=2\phi/2\pi\sim 0.64680~56887~25462~64456~47436~14345$.
The error is $\sim 10^{-31}$.  $\Q^{[1]}$ has not yet converged so a longer trajectory is needed.}
\label{ex2-conv}
\end{figure}

We study here the quasiperiodic orbit within a Siegel ball of the same system Eq. \ref{henon} but with different of parameter values $(\theta=0.664, \phi=2.032)$ and 
%starting from 
with an initial point $(x_0,y_0)=(-0.500+0.126i,-0.387-0.163i)\in \mathbb{C}^2$. 
Projections of the orbit are shown in Fig.~\ref{simple}, one of which is used to compute a rotation rate (Fig.~\ref{ex2-proj}~(Left). 
In order to compute the other rotation rate, unlike the case in Example 2A, here we use the ``time-2 delay'' $(r_n,r_{n-2})$ 
in Fig.~\ref{ex2-proj}~(Right).
Using these projections we measure angle differences, and the lifts of them are shown in Fig.~\ref{ex2-conv}~(a) and (b). 
The weighted Birkhoff averages ($\Q^{[1]}$ and $\Q^{[2]})$ of the lift of the angle difference (a) converges  to $\rho_1=(\theta+\phi)/2 \pi\sim 0.42908~17265~75749~82523~29106~26052$. 
%with an error of $10^{-31}.$ 
The lift of the angle difference (b) 
 converges  to $\rho_1-\rho_2=2\phi/2\pi\sim 0.64680~56887~25462~64456~47436~14345.$
%with an error of $10^{-32}.$ 
%

%\ifincludeXX
%\fi
%
%\ifincludeXX
%\fi
%
%\ifincludeXX
%\fi

\clearpage
\section{Computational methods}\label{s:method}
In the series of papers \cite{DCDS,EPL,DSSY,DY} we have developed techniques to characterize quasiperiodic orbits. 
We summarize them in this section used for our computation of this paper.
We introduced a weighted Birkhoff average ($\Q$) to calculate a Birkhoff average along a quasiperiodic orbit very quickly and in high accuracy \cite{EPL,DSSY}. 
The $\Q$ along a quasiperiodic orbit of the length $N$ is mathematically proved to show a faster convergence than any polynomials, when $N$ goes to infinity \cite{DSSY,DY}.

\subsection{The Birkhoff Ergodic Theorem}\label{ss:BET}
We use the Birkhoff Ergodic Theorem to prove the existence of a rotation vector 
and we use our weighted Birkhoff average to compute it and to compute Fourier coefficients of the torus $\torus$.
The Theorem assumes there is an invariant set, which here is the set $\torus$ with $d=1$ or $2$. Since we are interested here only in quasiperiodic dynamics, we can assume the dynamics are given by Eq.~\Xx\ where $\rho$ is irrational. 
Lebesgue measure is invariant; that is, each measurable set $E\subset \torus$ has the same measure as $F^{-1}(E) = E-\rho$ (and also the same as 
$F(E) = E+\rho$).
This map is ``ergodic'' because if $E$ is a set for which $E = F(E) = E+\rho$,
then the measure of $E$ is either $0$ or $1$. 

For an invariant measure $m$ enables the computation of the space-average $\int_\torus fdm$ for any $L^1$ function $f:\torus\to\R$ when a time series is the only information available. Since $m$ is Lebesgue measure, we can rewrite that integral as $\int_\torus f(\theta)dm$. 
We note that the Lebesgue measure of the entire torus is 1, so Lebesgue measure is a probability measure. Hence
$\int_\torus dm=1$.

For a map $F: \torus \to \torus$, the {\bf Birkhoff average} of
a function $f:\torus\to\R$ along the trajectory $\theta_n = F^n \theta_0$ is
\begin{equation}\label{eq:B}
B_N(f)(\theta_0) : = \frac{1}{N}\sum_{n=0}^{N-1} f(F^n(\theta_0))
\mbox{ and } L(f) := \int_X fdm.
\end{equation}
\allblack
\begin{teo}[Quasiperiodic case of the Birkhoff Ergodic Theorem \cite{BrinStuck}]\label{thm:ergodic_theorem} 

Let $F: \torus \to \torus$ satisfy Eq.~\ref {eq:plusrho} where $\rho\in\torus$ is irrational. Let $\mu$ be Lebesgue measure on $\torus$.  Then 
for every
%\protect\footnotemark
%\footnote
%\thefootnote
%\footnotetext{
%}
initial $\theta_0\in\torus$, 
$\lim_{N\to\infty}B_N(f)(\theta_0) =L(f)$.
\end{teo}
\allblue
The general ergodic theorem for general maps replaces the above ``for every'' with ``for almost every''. But for quasiperiodic maps the ``almost'' is unnecessary.
\allblack

For quasiperiodic dynamics the limit of $B_N$ as $N\to\infty$ exists for all $\theta_0$  and convergence is uniform in $\theta_0$ so we usually write $B_N(f)$ omitting $\theta_0$.
The Birkhoff average $B_N(f)$ 
%(\ref{eq:Birkhoff_sum})
 can be interpreted as an approximation to the integral $L(f)$ but convergence is very slow even for $C^\infty f$, in that for almost every rotation vector  there is a constant $C$ such that 
\[|B_N(f)- L(f)| \le C N^{-1},\] 
and even this slow rate will occur only under special circumstances such as when $( T^n(x))$ is a quasiperiodic trajectory. 
Let $f:\torus\to\R$ be $L^1$. 
\allblack

\subsection{The Method of Weighted Birkhoff Averages ($\Q^{[p]}_N$)}\label{ss:WB}
We have recently established a method for speeding up the convergence of 
the Birkhoff average in Theorem \ref{thm:ergodic_theorem} through introducing a $C^{\infty}$ weighting function
% by 
%\allred
%infinitely many
%\allblack
%orders of magnitude 
when the process is quasiperiodic and the function $f$ is $C^\infty$, a method we describe in ~\cite{EPL,DSSY,DY}. 
In \cite{DY} it is proved that the limit of using $\Q^{[p]}_N$ is the same as Birkhoff's limit.

Define the $C^{\infty}$ weighting function $w$ as follows.
\begin{equation}\label{eq:weight}
w^{[p]}(t) :=\begin{cases}
\exp\left(-[t(1-t)]^{-p}\right), & \mbox{for } t\in(0,1)\\
0, & \mbox{for } t\notin(0,1).
\end{cases}
\end{equation} 
Define the {\bf normalized weights} $\hat{w}^{[p]}_{n,N}:=\frac{w^{[p]}(n/N)}{\sum_{j=0}^{N-1}w^{[p]}(j/N)}$.
Note that $\sum_{n=0}^{N-1} \hat{w}^{[p]}_{n,N}=1$.
The {\bf Weighted Birkhoff  average} $\Q^{[p]}_N$ of $f$ is defined as follows.
\begin{equation}\label{eq:QN}
\Q^{[p]}_N(f)(\theta_0) :=\sum_{n=0}^{N-1} \hat{w}^{[p]}_{n,N}f(\theta_n).
\end{equation}

%In our calculations of the rotation rates, we use $p=1$ or $2$. 
% 
In our calculations of the rotation rates, we use $p=1$ or $2$. Sometimes this well-known function $w$ is referred to as a bump function because it is positive (with a single maximum) on $(0,1)$ and $w(x)$ and its derivatives of all orders are 0 at $x=0,1$.
See in particular \cite{DSSY} for details and a discussion of how the method relates to other approaches.
%\allred
We remark that we generally observe faster convergence using $p=1$ for up to 10-digit accuracy but higher precision $p=2$ converges faster.
In even higher precision computation (greater than 40 digits) we can get faster convergence by applying the normalized weighting function of 
\begin{equation}\label{eq:weighthigh}
w^{[\ast]}(t) :=\begin{cases}
\exp\left(-1/w^{[1]}\right), & \mbox{for } t\in(0,1)\\
0, & \mbox{for } t\notin(0,1).
\end{cases}
\end{equation} 
%\allblack

\subsection{A ``super convergence'' theorem for $\Q$.} 
\begin{teo}\label{thm:A1}
Let $X$ be a $C^\infty$ manifold and $T:X \to X$ be a d-dimensional $C^\infty$ quasiperiodic map on $X_0\subseteq X$, with invariant probability measure $m$. 
Assume $T$ has a Diophantine rotation vector.
Let $f:X \to\R$ be $C^\infty$.
Assume $w^{[p]}$ is the 
%exponential 
weighting function (in Eq.~(\ref{eq:weight})) where $p\ge 1$. 
%\allred
Write $L^{(N)}(f):=\Q_N^{[p]}(f)(x_0).$
Then for each $x_0\in X_0$, 
\begin{equation}
%\mbox{ For each }m \mbox{ there exists }N_m >0 \mbox{ such that } |\Q_N^{[p]}(f)(x_0) - L(f)| \le N^{-m} \mbox{ for all }N>N_m.\label{super}
 \mbox{ for each }  m\in\mm N,~
 |L^{(N)}(f) - L(f)|~N^{m} \to 0 \mbox{ as } N\to\infty.\label{super}
\end{equation}
%\allblack
\end{teo}
We call the above Property \ref{super} ``super convergence''. 
See \cite{DY} for a more general statement.
See in particular \cite{DSSY} for details and a discussion of how the method relates to other approaches.

Essentially the same weighting function  as $w^{[1]}$ is discussed by Laskar~\cite{Laskar99} in the Remark~2 of the Annex (i.e. Appendix) without proofs
and without numerical examples.

In \cite{Laskar99,Laskar93a,Laskar93b,Laskar03}, Laskar computes using a Hanning data weighting function, essentially using the weighting function  $w_{\sin^2}=\sin^2$  for his computations, which leads to the  convergence rate of order approximately $N^{-2.5}$ (where $N$ is the length of a orbit), as observed in \cite{DSSY}. See that paper's discussion and Fig. 7.
Our Fourier series calculation required a day on one processor. Had we used his $w_{\sin^2}$ method, it would have required 1 million days.

See also the discussion in  \cite{DSSY} of the methods of other authors, particularly \cite{seara:villanueva:06, luque:Rot, luque:villanueva:14}. These methods are fast but require a complicated iterative scheme.

\subsection{Computing Fourier coefficients using $\Q$.}
% the weighted Birkhoff average $\Q$.}
Let $(\theta_n)$ be a quasiperiodic trajectory on $[0,1)$.
Each Fourier coefficient of a function $g:\mathbb{C}\to\mathbb{C}$ can be written $b_k:=\displaystyle\int_0^1  g(e^{i2\pi\theta})\sigma_{-k}(\theta)d\theta$. 
We compute $\hat{b}_k$ $:=\Q_N(g(e^{i2\pi\theta})\sigma_{-k}(\theta))$ to approximate $b_k$; 
substituting  $\displaystyle\sum_j b_j \sigma_j$ for $ g(e^{i2\pi\theta})$ gives, 
$$
\hat{b}_k=\displaystyle\sum_j~b_j \Q_N(\sigma_{j-k})=\displaystyle\sum_j~b_{j+k} \Q_N(\sigma_{j}).
$$
The trajectory average of $\sigma_{j}$ is $0$ except for $j=0$. 
When 
\allblack
$\left|\Q_N(\sigma_j)\right|<\varepsilon$, 
the error in estimate is
$$\left|\hat{b}_k-b_k\right|=\ \left|\displaystyle\sum_{j\neq 0}~b_{j+k} \Q_N(\sigma_{j})\right|\ \approx \varepsilon.$$
\section{Summary}\label{s:summary}
We have investigated quasiperiodic orbits in one-dimensional and two-dimensional complex dynamical systems. 
The usefulness of the weighted average and the many possibilities for rotation vectors in a high dimensional quasiperiodic orbit 
are especially exemplified. 
\section*{Acknowledgments}
The authors would like to thank Prof. Shigehiro Ushiki for providing us much useful information about the Siegel ball of the complex H\'enon map. 

%\red{I MOVED THE SUMMATION APPENDIX TO AFTER THE END OF DOCUMENT.}

\bibliographystyle{unsrt}%{abbrv}{iopart-num}
\bibliography{qr_bibliography,Weighted_calc_bibliography}
\end{document}